\documentclass[11pt]{article}
\usepackage{float}
\usepackage{url}
\usepackage{mathrsfs}
\usepackage{color}
\usepackage[toc,page]{appendix} 

\usepackage{latexsym}
\usepackage{amsmath}
\usepackage{amsfonts}
\usepackage{amssymb}
\usepackage{amsthm}
\usepackage{psfrag}
\usepackage{epsfig}
\usepackage[english]{babel}
\usepackage{color}
\usepackage{soul}
\usepackage{marginnote}

\usepackage{tikz}
\usepackage{pgfplots}
\usepackage{pgfplotstable}
\usetikzlibrary{calc}

\newcommand{\slopeTriangleBelow}[6]
{
    \pgfplotsextra{
    \coordinate (A) at (axis cs:#1,#3);
    \coordinate (B) at (axis cs:#1,#4);
    \coordinate (C) at (axis cs:#2,#4);

    \draw[#6] (A)--(B) node [midway, left] {#5} ;
    \draw[#6] (B)--(C)  node [midway, above] {1};
    \draw[#6] (C)--(A) ;
    }
}

\newcommand{\slopeTriangleAbove}[6]
{
    \pgfplotsextra{ 
    \coordinate (A) at (axis cs: #1,#3);
    \coordinate (B) at (axis cs: #2,#3);
    \coordinate (C) at (axis cs: #2,#4);

    \draw[#6] (A)--(B) node [midway, below]{1};
    \draw[#6] (B)--(C) node [midway, right]{#5};
    \draw[#6] (C)--(A) ;
    }
}

\newsavebox{\fmbox}

\newtheorem{theorem}{Theorem}[section]
\newtheorem{lemma}[theorem]{Lemma}

\newtheorem{remark}[theorem]{Remark}

\newtheorem{proposition}[theorem]{Proposition}

\newtheorem*{remark*}{Remark}
\newtheorem*{assumption*}{Assumption}

\newcommand{\pf}{{\bf Proof: }}

\newcommand{\eps}{\varepsilon}

\newcommand{\dis}{\displaystyle}

\def\Vo{\vbox{\offinterlineskip\hbox{\kern 3pt$\scriptstyle\circ$}
\kern 1pt\hbox{$V$}}}
\def\Ho{\vbox{\offinterlineskip\hbox{\kern 3pt$\scriptstyle\circ$}
\kern 1pt\hbox{$H$}}}
\def\Wo{\vbox{\offinterlineskip\hbox{\kern 3pt$\scriptstyle\circ$}
\kern 1pt\hbox{$W$}}}

\newcommand{\R}{\mathbb{R}}

\newcommand{\caC}{{\cal C}}
\newcommand{\caD}{{\cal D}}

\newcommand{\caK}{{\cal K}}
\newcommand{\caA}{{\cal A}}
\newcommand{\caQ}{{\cal Q}}

\renewcommand{\textbf}[1]{{\bfseries\mathversion{bold}#1}}
\newcommand{\bfu}{\boldsymbol{u}}

\newcommand{\s}{{\sigma}}
\newcommand{\sn}{{\sigma_n}}

\newcommand{\vn}{v_n}

\renewcommand{\div} {{\rm{div} \,}}
\newcommand{\GC}{\Gamma_C}

\newcommand{\GCQ}{Q_C}
\newcommand{\GCQtilde}{\tilde{Q}_C}

\newcommand{\ZC}{{Z_{C}(\GCQ)}}
\newcommand{\ZNC}{{Z_{NC}(\GCQ)}}

\newcommand{\intGCQ}{\int_{\GCQ}}

\newcommand{\dG}{\ {\rm d} \Gamma}
\newcommand{\dO}{\ {\rm d}\Omega}
\newcommand{\intGC}{\int_{\Gamma_C}}
\newcommand{\intGN}{\int_{\Gamma_N}}
\newcommand{\intO}{\int_{\Omega}}

\newcommand{\abs}[1]{\left|{#1}\right|}
\newcommand{\norm}[1]{\left\|{#1}\right\|}

\newcommand{\scal}[2]{\left<#1,#2\right>_{W',W}}

\newcommand\restr[2]{{% we make the whole thing an ordinary symbol
  \left.\kern-\nulldelimiterspace % automatically resize the bar with \right
  #1 % the function
  \vphantom{\big|} % pretend it's a little taller at normal size
  \right|_{#2} % this is the delimiter
  }}
\newcommand\restrr[2]{\ensuremath{\left.#1\right|_{#2}}}
\newcommand{\cqfd}{{$\mbox{}$\hfill$\square$}}

\usepackage{letltxmacro}
\makeatletter
\let\oldr@@t\r@@t
\def\r@@t#1#2{%
\setbox0=\hbox{$\oldr@@t#1{#2\,}$}\dimen0=\ht0
\advance\dimen0-0.2\ht0
\setbox2=\hbox{\vrule height\ht0 depth -\dimen0}%
{\box0\lower0.4pt\box2}}
\LetLtxMacro{\oldsqrt}{\sqrt}
\renewcommand*{\sqrt}[2][\ ]{\oldsqrt[#1]{#2} }
\makeatother

\renewcommand{\textbf}[1]{{\bfseries\mathversion{bold}#1}}
\newcommand{\bfz}{\boldsymbol{\zeta}}
\newcommand{\bfX}{\boldsymbol{\Xi}}

\newcommand{\bfi}{\boldsymbol{i}}

\newcommand{\bfI}{\boldsymbol{I}}

\newcommand{\bflambda}{\boldsymbol{\lambda}}
\newcommand{\bfN}{\boldsymbol{N}}
\newcommand{\bfB}{\boldsymbol{B}}

\usepackage{subfigure}

\usepackage[normalem]{ ulem }
 \usepackage{soul}

\title{\textbf{ \textit{A priori} error for unilateral contact problems with Lagrange multipliers and IsoGeometric Analysis }}
\author{Pablo Antolin\footnote{EPFL SB MATHICSE MNS (B\^at. MA), station 8, CH 1015 Lausanne (Switzerland).}~,
Annalisa Buffa$^{*,}$\footnote{\samepage Istituto di Matematica Applicata e Tecnologie Informatiche 'E. Magenes' del CNR
via Ferrata 1, 27100, Pavia (Italy). \qquad\qquad\qquad\qquad\qquad\qquad\qquad\qquad\qquad\qquad\qquad\qquad\qquad\qquad\qquad\qquad\qquad\qquad\qquad\qquad\qquad\qquad $ $  email: pablo.antolin@epfl.ch,  annalisa.buffa@epfl.ch,  mathieu.fabre@epfl.ch.}~,
Mathieu Fabre$^{*,\dag}$}

\begin{document}

\renewcommand{\labelitemi}{\textbullet}  
\maketitle
%\tableofcontents
 
 \begin{abstract}
In this paper, we consider unilateral contact problem without friction between a rigid body and deformable one in the framework of isogeometric analysis. We present the theoretical analysis of the mixed problem{. For the displacement, we use the pushforward of a NURBS space of degree $p$ and for the Lagrange multiplier, the pushforward of a B-Spline space of degree $p-2$. These chooses of space ensure to prove an $\inf-\sup$ condition and so on, the stability of the method. An active set strategy is used in order to avoid of geometrical hypothesis of the contact set.}
 An {optimal} \textit{a priori} error estimate is demonstrated without assumption on the unknown contact set. Several numerical examples in two- and three-dimensions and in small and large deformation demonstrate the accuracy of the proposed method. 
\end{abstract}
 
 %List of figures 
 %\listoffigures
 
%%%%%%%%%%%%%%%%%%%%%%%%%%%%%%%%%%%%%%%%%%%%%%%% 
 \section*{Introduction}
\label{sec:sec0}
\addcontentsline{toc}{section}{Introduction}
 
In the past few years, the study of contact problems in small and large deformation is increased. The numerical solution of contact problems presents several difficulties as the computational cost, the high nonlinearity and the ill-conditioning. Contrary to many other problems in nonlinear mechanics, these problems can not be solved always at a satisfactory level of robustness and accuracy \cite{laursen-02,wriggers-06} with the introduce numerical methods.

One of the reasons that make robustness and accuracy hard to achieve is that the computation of gap, i.e. the distance between the deformed body and the obstacle is indeed an ill-posed problem and its numerical approximation often introduce extra discontinuity that breaks the converge of the iterative schemes; see \cite{al-cu1988,laursen-02,wriggers-06,kon-sch-13} where a master-slave method is introduced to weaken this effect.

To this respect, the use of NURBS or spline approximations within the framework of isogeometric analysis \cite{hughes05}, holds great promises thanks to the increased regularity in the geometric description which makes the gap computation intrinsically easier. 
{The IGA-based methods use a generalization of B\'ezier's curves, the B-Splines and non-uniform rational B-Splines (NURBS). These functions, used to represent the geometry of the domains with CAD, are used as basis functions to approximate a partial differential equations, it is called the isoparametric paradigm. The smooth IGA basis functions possess a number of signifiant advantages for the analysis, including exact geometry and superior approximation. }
Isogeometric methods for \hbox{frictionless} contact problems have been introduced in \cite{lorenzis9,Temizer11,Temizer12,lorenzis12,lorenzisreview,lorenzis15}, see also with {primal} and dual elements \cite{Wohlmuth00,hueber-Wohl-05,hueber-Stadler-Wohlmuth-08,wohl12,seitz16}. Both point-to-segment and segment-to-segment (i.e, mortar type) algorithms have been designed and tested with an engineering perspective, showing that, indeed, the use of smooth geometric representation helps the design of reliable methods for contact problems.

In this paper, we take a slightly different point of view. Inspired by the recent design and analysis of isogeometric mortar methods in \cite{brivadis15}, we consider a formulation of frictionless contact based on the choice of the Lagrange multiplier space proposed there. Indeed, we associate to NURBS displacement of degree $p$, a space of Lagrange multiplier of degree $p-2$. The use of lower order multipliers has several advantages because it makes the evaluation of averaged gap values at active and inactive control points simpler, accurate and substantially more local.  This choice of multipliers is then coupled with an active-set strategy, as the one proposed and used in \cite{hueber-Wohl-05,hueber-Stadler-Wohlmuth-08}. 

Finally, we perform a comprehensive set of tests both in small and large scale deformation, which well show the performance of our method. These tests have been performed with an in-house code developed upon the public library igatools \cite{Pauletti2015}.

The outline of the paper is structured as follows in Section \ref{sec:sec1}, we introduce unilateral contact problem, some notations. In Section \ref{sec:discrete_pb}, we describes the discrete spaces and their properties.
% and the discrete spaces%, NURBS of degree $2$ for the displacement and B-Spline of degree $0$ for the Lagrange Multiplier
In Section \ref{sec:sec2}, we present the theoretical analysis of the mixed problem. An {optimal} \textit{a priori} error estimate without assumption on the unknown contact set is presented.
%we define the Lagrange multiplier method together with active-set strategy and some numerical results for NURBS of degree $2$ for the displacement and B-Spline of degree $0$ for the Lagrange Multiplier. %The case of NURBS of degree $3$ for the displacement and B-Spline of degree $1$ for the Lagrange Multiplier does not converge with this method.
In the last section, some two- and three-dimensional problem in small deformation are presented in order to illustrate the convergence of the method with active-set strategy. A two-dimensional problem in large deformation with Neo-Hookean material law is provided to show the robustness of this method.

\begin{remark*}
The letter $C$ stands for a generic constant, independent of the discretization parameters and the solution $u$ of the variational problem. For two scalar quantities $a$ and $b$, the notation $a \lesssim b$ means there exists a constant $C$, independent of the mesh size parameters, such that $a \leq Cb$. Moreover, $a \sim b$ means that $a \lesssim b$ and $b \lesssim a$.
\end{remark*}

%%%%%%%%%%%%%%%%%%%%%%%%%%%%%%%%%%%%%%%%%%%%%%%% 
  
\section{Preliminaries and notations}
\label{sec:sec1}
\subsection{Unilateral contact problem}
\label{subsec:continuous_pb}

Let $\Omega \subset \R^d$ ($d = $2$ \textrm{ or } 3$) be a bounded regular domain which represents the reference configuration of a body submitted to a Dirichlet condition on $\Gamma_D$ (with $\textrm{meas}(\Gamma_D) >0$), a Neumann condition on $\Gamma_N$ and a unilateral contact condition on a potential zone of contact $\Gamma_C$ with a rigid body. Without loss of generality, it is assumed that the body is subjected to a volume force $f$, to a surface traction $\ell$ on $\Gamma_N$ and clamped at $\Gamma_D$. Finally, we denote by $n_{\Omega}$ the unit outward normal vector on $\partial \Omega$.

In what follows, we call $u$ the displacement of $\Omega$, $\dis \eps (u) = \frac{1}{2} (\nabla u + \nabla u^T)$ its linearized strain tensor and we denote by $\dis \s  = (\s_{ij})_{1 \leq i,j \leq d} $ the stress tensor. We assume a linear constitutive law between $\s$ and $\eps$, \textit{i.e.}  $\dis \s (u)= A \eps (u)$, where $A=(a_{ijkl})_{1 \leq i,j,k,l \leq d}$ is a fourth order symmetric tensor verifying the usual bounds:
\begin{itemize}
\item $\dis a_{ijkl} \in L^{\infty}(\Omega), \textrm{ \textit{i.e.} there exists a constant $m$ such that } \max_{1 \leq i,j,k,l \leq d} \abs{a_{ijkl}} \leq m;$
\item $\dis \textrm{there exists a constant } M>0 \textrm{ such that \textit{a.e.} on } \Omega,$
$$\dis a_{ijkl}  \eps_{ij}\eps_{kl}  \geq M  \eps_{ij}   \eps_{ij}  \quad \forall \eps \in \R^{d \times d} \textrm{ with } \eps_{ij} = \eps_{ji}.$$
\end{itemize} 

   Let $n$ be the outward unit normal vector at the rigid body. {From now on, we assume that $n$ is an infinitely regular field}.  For any 
 %The linearized strain tensor field is given by $\dis \eps (v) = \frac{1}{2} (\nabla v + \nabla v^T)$, the stress tensor field $\dis \s  = (\s_{ij})_{1 \leq i,j \leq d} $, $A$ is the fourth order symmetric elasticity tensor satisfying the usual uniform ellipticity and boundedness properties , $\tilde{n}$ is the outward unit normal vector of the body. 
%and ${n}$ is the outward unit normal vector of the rigid body. For any 
displacement field $u$ and for any density of surface forces $\s (u) n$ defined on $\partial \Omega$, we adopt the following notation: $$ u = u_n n + u_t \qquad \textrm{and} \qquad \s (u) n = \s_n (u) n + \s_t (u) ,$$ where $u_t$ (resp. $\s_t (u)$) are the tangential components with respect to $n$.% the outward unit normal vector of the rigid body, of $u$ (resp. $\s (u) n$). 

The unilateral contact problem between a rigid body and the elastic body $\Omega$ consists in finding the displacement $u$ satisfying:
\begin{eqnarray} \label{eq:strong}
\begin{array}{rcll}
 \div \s (u) +f \!\!\!&\!\!\!=\!\!\!&\!\!\! 0  &\qquad \textrm{in } \Omega, \\
 \s (u) \!\!\!&=&\!\!\! A \eps (u)&\qquad \textrm{in } \Omega,\\
 u \!\!\!&=&\!\!\! 0 &\qquad \textrm{on } \Gamma_D,\\
  \s (u) {n_\Omega} \!\!\!&=&\!\!\! \ell &\qquad \textrm{on } \Gamma_N.
 \end{array}
\end{eqnarray} 
and the conditions describing unilateral contact without friction at $\Gamma_C$ are:
\begin{eqnarray} \label{eq:contact_cond}
\begin{array}{rl}
 u_n \geq &\!\!\! 0  \quad (i), \\
 \s_n (u) \leq &\!\!\! 0  \quad (ii), \\
  \s_n (u) u_n = &\!\!\! 0  \quad (iii), \\
 \s_t (u) = &\!\!\! 0  \quad (iv).
\end{array}
\end{eqnarray} 

\noindent In order to describe the variational formulation of \eqref{eq:strong}-\eqref{eq:contact_cond}, we consider the Hilbert spaces: $$ \dis V := H^1_{0,\Gamma_D}(\Omega)^d = \{ v \in H^1(\Omega)^d , \quad v = 0 \textrm{ on } \Gamma_D \}, \quad W=\{ \restrr{v_n}{\Gamma_C}, \quad v \in V \},$$ and their dual spaces $V'$, $W'$ endowed with their usual norms. We denote by: $$\dis \norm{v}_{V} = \left(\norm{v}_{L^2(\Omega)^d}^2 + \abs{v}_{H^1(\Omega)^d}^2\right)^{1/2} , \ \forall v \in V.$$
If $\overline{\Gamma}_D \cap \overline{\Gamma}_C = \emptyset$ and $n$ is regular enough, it is well known that $W=H^{1/2}(\Gamma_C) $ and we denote $W'$ by $H^{-1/2}(\Gamma_C)$. On the other hand, if $\overline{\Gamma}_D \cap \overline{\Gamma}_C \neq \emptyset$, it will hold that $H^{1/2}_{00} (\Gamma_C) \subset W \subset H^{1/2}(\Gamma_C)$.

\noindent In all cases, we will denote by $\norm{\cdot}_W $ the norm on $W$ and by $\scal{\cdot}{\cdot}$ the duality pairing between $W'$ and $W$.
%and $ \scal{.}{.}$ represents the duality pairing between $H^{-1/2}(\Gamma_C)$ and $H^{1/2}(\Gamma_C)$. As a reminder, we have $H^{1/2}_{00}(\Gamma_C) \subset W \subset H^{1/2}(\Gamma_C)$ and $W' \subset H^{-1/2}(\Gamma_C)$. Classically, $ H^{1/2}(\Gamma_C)$ is the space of restrictions on $\Gamma_C$ of traces on $\partial \Omega$ of functions in $H^1(\Omega)^d$ and $H^{-1/2}(\Gamma_C)$ is the dual space of $H^{1/2}_{00}(\Gamma_C)$ which is the space of the restrictions on $\Gamma_C$ of functions in $H^{1/2}(\Gamma_C)$ vanishing outside $\Gamma_C$. For more details, we can refer to \cite{lions-magenes-72}.

\noindent For all $u$ and $v$ in $V$, we set: $$\dis a(u,v) = \intO \s(u): \eps(v)\dO \quad \textrm{and} \quad L(v) = \intO f \cdot v \dO + \intGN \ell \cdot v \dG.$$

\noindent Let ${K_C}$ be the closed convex cone of admissible displacement fields satisfying the non-interpenetration conditions, $\dis {K_C} := \{ v \in V, \quad v_n \geq 0 \textrm{ on } \Gamma_C\}$. {A weak formulation of Problem} \eqref{eq:strong}-\eqref{eq:contact_cond} (see \cite{lions-magenes-72}), as a variational inequality, is to find $u \in {K_C} $ such as: 
\begin{eqnarray} \label{eq:var_ineq}
\dis a(u,v-u) \geq L(v-u), \qquad \forall v \in {K_C}.
\end{eqnarray} 

We cannot directly use a Newton-Raphson's method to solve the formulation \eqref{eq:var_ineq}. A classical solution is to introduce a new variable, the Lagrange multipliers denoted by $\lambda$, which represents {the surface normal force} % { or which takes care of the non-penetration constraints}
. For all $\lambda$ in $W'$, we denote $\dis b(\lambda,v) = - \scal \lambda \vn$ and $M$ is the classical convex cone of multipliers on $\Gamma_C$:
 $$ \dis M := \{ \mu \in W', \quad \scal{\mu}{\psi} \leq 0 \quad \forall \psi \in H^{1/2}(\Gamma_C), \quad \psi \geq 0 \textit{ a.e.} \textrm{ on } \Gamma_C \} .$$
 \noindent The complementary conditions with Lagrange multipliers writes as follows:
\begin{eqnarray} \label{eq:contact_cond_mult}
\begin{array}{rl}
 u_n \geq &\!\!\! 0  \quad (i), \\
\lambda \leq &\!\!\! 0  \quad (ii), \\
 \lambda u_n = &\!\!\! 0  \quad (iii).
\end{array}
\end{eqnarray} 

\noindent {The mixed formulation \cite{ben-belgacem-renard-03} of the Signorini problem \eqref{eq:strong} and \eqref{eq:contact_cond_mult} consists in finding $(u,\lambda) \in V \times M$ such that:}
\begin{eqnarray} \label{eq:mixed_form}
\left\{
\begin{array}{rl}
\dis a(u,v) - b(\lambda,v) = L(v),& \dis \qquad \forall v \in V,\\
\dis b(\mu-\lambda,u) \geq 0,&\dis \qquad \forall \mu \in M.
\end{array}
\right.
\end{eqnarray} 
Stampacchia's Theorem ensures that problem \eqref{eq:mixed_form} admits a unique solution.
  
%We can too resolve a saddle-point as follows finding $(u,\lambda) \in V \times M$ such that: 
%\begin{eqnarray} \label{eq:saddle_point}
%\dis \caL(u,\mu) \leq \caL(u,\lambda) \leq \caL(v,\lambda) , \qquad \forall v \in V , \forall \mu \in M ,
%\end{eqnarray} 
%where $ \caL(.,.)$ is the classical Lagrangian such that $ \caL(v,\mu) = \frac 1 2 a(v,v) - L(v) - b(\mu,v)$.
\noindent The existence and uniqueness of the solution $(u,\lambda)$ of the mixed formulation has been established in \cite{haslinger-96} and it holds $\lambda = \sn (u)$. 

\noindent {To simplify the notation, we denote by $\norm{\cdot}_{3/2+s,\Omega} $ the norm on $H^{3/2+s}(\Omega)^d$ and by $\norm{\cdot}_{s,\GC} $ the norm on $H^s(\Gamma_C)$.}

\noindent So, the following classical inequality (see \cite{daveiga06}) holds:
\begin{theorem}\label{thm:u-lambda}
Given $s>0$, if the displacement $u$ verifies $u \in H^{3/2+s}(\Omega)^d$, then $\lambda \in H^s(\Gamma_C)$ and it holds:
\begin{eqnarray} \label{ineq:Sp-2 primo}
 \dis \norm{\lambda}_{s,\GC} \leq  \norm{u}_{3/2+s,\Omega}.
\end{eqnarray} 
\end{theorem}
%To approximate the mixed formulation \eqref{eq:mixed_form}, we use a work presented in \cite{brivadis15} on the mortar method. A $\inf-\sup$ stability is proved using a primal space of NURNS of degree $p$ and dual space of B-Spline of degree $p-2$. On the next of the purpose, we use these spaces as approximation spaces.
The aim of this paper is to discretize the problem \eqref{eq:mixed_form} within the isogeometric paradigm, \textit{i.e.} with splines and NURBS. Moreover, in order to properly choose the space of Lagrange multipliers, we will be inspired by \cite{brivadis15}. In what follows, we introduce NURBS spaces and assumptions together with relevant choices of space pairings. In particular, following \cite{brivadis15}, we concentrate on the definitions of B-Splines displacements of degree $p$ and multiplier spaces of degree $p-2$.

\subsection{NURBS discretisation}
\label{subsec:nurbs_discretisation}

In this section, we describe briefly an overview on isogeometric analysis providing the notation and concept needed in the next sections. Firstly, we define B-Splines and NURBS in one-dimension. Secondly, we extend these definitions to the multi-dimensional case. Finally, we define the primal and the dual spaces for the contact boundary.

{Let us denote by $p$ the degree of univariate B-Splines and by $\Xi$ an open univariate knot vector, where the first and last entries are repeated $(p + 1)$-times, \textit{i.e.} $$\Xi := \{ 0= \xi_1= \cdots =  \xi_{p+1} <  \xi_{p+2} \leq \ldots \leq  \xi_{\eta}< \xi_{\eta+1}= \cdots =  \xi_{\eta+p+1} \}. $$}
Let us define $Z = \{ \zeta_1, \ldots, \zeta_E\}$ as vector of breakpoints, \textit{i.e.} knots taken without repetition, and $m_j$, the multiplicity of the breakpoint $ \xi_j, \ j=1, \ldots , E$.  Let $\Xi$ be the open knot vector associated to $Z$ where each breakpoint is repeated $m_j$-times, \textit{i.e.} 
In what follows, we suppose that $m_1 = m_E = p+1$, while $m_j \leq p-1$, $\forall j = 2,\ldots, E-1$. We define by $\hat{B}^p_i(\zeta)$, $i=1,\ldots,\eta$ the $i$-th univariable B-Spline based on the univariate knot vector $\Xi$ and the degree $p$. We denote by $S^p (\Xi)=  Span\{ \hat{B}^p_i(\zeta), \ i=1,\ldots,\eta \}$. Moreover, for further use we denote by $\tilde{\Xi}$ the sub-vector of $\Xi$ obtained by removing the first and the last knots.
% the corresponding B-Spline space and $\Lambda^h (\tilde{\Xi})=  Span\{ \hat{B}^p_{\lambda,i}(\zeta), \ i=1,\ldots,\eta \}$. It is important to underline that each basis function is $C^{p-m_j}$ at $\zeta_j$, $j=1,\ldots,E$.
%
%\noindent Let us denote by $\{\omega_i\}_i$ as a given set of positive weights and $\dis \hat{W}(\zeta) = \sum_{i=1}^n \omega_i \hat{B}^p_i(\zeta) $ as weights function. We define the NURBS (Non Uniform Rational Basis Splines) as $$\dis \hat{N}^p_i(\zeta) = \frac{\omega_i  \hat{B}^p_i(\zeta)}{\hat{W}(\zeta)}, $$ and let $N^p(\Xi)$ be as the NURBS space.\\

Multivariate B-Splines in dimension $d$ are obtained by tensor product of univariate B-Splines. For any direction $\delta \in \{ 1, \ldots, d\}$, we define by 
%$p_\delta$ the degree of the B-Splines, 
$\eta_\delta$ the number of B-Splines, $\Xi_\delta$ the open knot vector and $Z_\delta$ the breakpoint vector. Then, we define the multivariate knot vector by $\bfX = ( \Xi_1 \times \ldots \times \Xi_d) $ and the multivariate breakpoint vector by $\boldsymbol{Z} = ( Z_1 \times \ldots \times Z_d ) $. We introduce a set of multi-indices $\bfI = \{ \bfi  =(i_1, \ldots, i_d) \mid 1 \leq i_\delta \leq \eta_\delta \}$. %To simplify the notation, we take $p_\delta = p$. 
We build the multivariate B-Spline functions for each multi-index $\bfi$ by tensorization from the univariate B-Splines, {let $\bfz \in \boldsymbol{Z}$ be a parametric coordinate of the generic point}: $$\hat{B}^p_{\bfi}(\bfz) = \hat{B}^p_{i_1}(\zeta_1) \ldots \hat{B}^p_{i_d}(\zeta_d) .$$
Let us define the multivariate spline space in the reference domain by (for more details, see \cite{brivadis15}): $$S^p(\bfX) =  Span\{ \hat{B}^p_{\bfi}(\bfz), \ \bfi\in \bfI \}.$$
We define $N^p(\bfX)$ as the NURBS space, spanned by the function $\hat{N}^p_{\bfi}(\bfz)$ with $$\hat{N}^p_{\bfi}(\bfz)  = \frac{\omega_{\bfi}  \hat{B}^p_{\bfi}(\bfz)}{\hat{W}(\bfz)},  $$ {where $\{\omega_{\bfi}\}_{\bfi\in \bfI}$ is a set of positive weights and $\dis \hat{W}(\bfz) = \sum_{\bfi\in \bfI} \omega_{\bfi} \hat{B}^p_{\bfi}(\bfz) $ is the weight function and we set }$$ N^p(\bfX) =  Span\{ \hat{N}^p_{\bfi}(\bfz) , \ \bfi\in \bfI  \} .$$

In what follows, we will assume that $\Omega$ is obtained as image of $\dis\hat{\Omega} = ]0,1[^d$ through a NURBS mapping $\varphi_0$, \textit{i.e.} $\dis\Omega = \varphi_0(\hat{\Omega})$. Moreover, in order to simplify our presentation, we assume that $\Gamma_C$ is the image of a full face $\dis\hat{f}$ of $\dis\bar{\hat{\Omega}}$, \textit{i.e.} $\dis{\Gamma_C} = \varphi_0(\hat{f})$. We denote by $\dis \varphi_{0,\Gamma_C}$ the restriction of $\dis \varphi_{0}$ to $\dis\hat{f}$. \\

{
\noindent  A NURBS surface, in d=2, or solid, in d=3, is parameterised by $$ \caC( \bfz) = \sum_{\bfi \in \bfI} C_{\bfi} \hat{N}^p_{\bfi}(\bfz) ,$$
where ${C_{\bfi}}_{ \in \bfI} \in \R^d$, is a set of control point coordinates. The control points are somewhat analogous to nodal points in finite element analysis. The NURBS geometry is defined as the image of the reference domain $\hat{\Omega}$ by $\varphi$, called geometric mapping, $\Omega_t = \varphi(\hat{\Omega})$. \\
}

We remark that the physical domain $\Omega$ is split into elements by the image of $\boldsymbol{Z}$ through the map $\varphi_0$. We denote such a {physical mesh $\caQ_h$ }and {physical elements} in this mesh will be called $Q$. %\noindent For any $Q \in \caQ_h$, $\tilde{Q}$ denotes the support extension of Q (see \cite{daveiga06,daveiga2014}) defined as the image of supports of B-Splines that are not zero on $\hat{Q}= \varphi_0^{-1} (Q)$. 
$\Gamma_C$ inherits a mesh that we denote by $ \restr{\caQ_h}{\GC} $. Elements on this mesh will be denoted as $Q_C$. \\

Finally, we introduce some notations and assumptions on the mesh. 

\noindent \textbf{Assumption 1.} The mapping $\varphi_0$ is considered to be a bi-Lipschitz homeomorphism. Furthermore, for any parametric element ${\hat{Q}}$, $ \dis \restr{\varphi_0}{\bar{\hat{Q}}}$ is in $\caC^\infty(\bar{\hat{Q}})$ and for any physical element ${{Q}}$, $ \dis \restr{\varphi_0^{-1}}{{\bar{Q}}}$ is in $\caC^\infty({\bar{Q}})$.

\noindent Let $h_Q$ be the size of an {physical element} $Q$, it holds $h_Q = \textrm{diam} (Q)$. In the same way, we define the mesh size for any parametric element. In addition, the Assumption 1 ensures that both size of mesh are equivalent. We denote the maximal mesh size by $\dis h = \max_{Q \in \caQ_h} h_Q$.

\noindent \textbf{Assumption 2.} The mesh $\caQ_h$ is quasi-uniform, \textit{i.e} there exists a constant $\theta$ such that $\dis \frac{h_Q}{h_{Q'}} \leq \theta$ with $Q$ and $Q' \in \caQ_h$.

\section{Discrete spaces and their properties}
\label{sec:discrete_pb}
%To approximate the problem \eqref{eq:mixed_form}, first we define the discretized spaces to define then the formulation.\\
We concentrate now on the definition of spaces on the domain $\Omega$. 

\noindent For displacements, we denote by $V^h \subset V$ the space of mapped NURBS of degree $p$ with appropriate homogeneous Dirichlet boundary condition: $$\dis V^h := \{ v^h = \hat{v}^h \circ \varphi^{-1}_0,  \quad \hat{v}^h \in N^p(\bfX)^d \} \cap V .$$
%\noindent Let us introduce the approximations spaces. Let $V^h \subset V$ be: $$\dis V^h = \{ v^h = \hat{v}^h \circ \varphi^{-1}_0, \quad \hat{v}^h \in N^p(\bfX) \} \cap V ,$$ %$V^h = \{ v^h = \hat{v}^h \circ \varphi^{-1}_0, \quad \hat{v}^h \in N^p(\bfX), \quad v^h = 0 \textrm { on } \GD \}$ 
%defined on the knot vector $\bfX$ of degree $p$. We denote by $h$ the mesh size of $V^h$.

\noindent We denote the space of traces normal to the rigid body as: $$\dis W^h := \{ \psi^h, \quad \exists v^h \in V^h: \quad v^h \cdot n = \psi^h \textrm{ on } \GC  \}.$$

\noindent For multipliers, following the ideas of \cite{brivadis15}, we define the space of B-Splines of degree $p-2$ on the potential contact zone $\dis \Gamma_C = \varphi_{0,\Gamma_C}(\hat{f})$. We denote by $\bfX_{\hat{f}}$ the knot vector defined on $\hat{f}$ and by $\tilde{\bfX}_{\hat{f}}$ the knot vector obtained by removing the first and last value in each knot vector. We define:
$$\dis \Lambda^h := \{ \lambda^h = \hat{\lambda}^h \circ \varphi_{0,\Gamma_C}^{-1}, \quad \hat{\lambda}^h \in S^{p-2}(\tilde{\bfX}_{\hat{f}}) \} .$$
The {scalar} space $\Lambda^h$ is spanned by mapped B-Splines of the type $\dis \hat{B}^{p-2}_{\bfi}(\bfz)\circ \varphi_{0,\Gamma_C}^{-1}$ for $\bfi$ belonging to a suitable set of indices. In order to reduce our notation, we call {$K$ the unrolling of the multi-index $\bfi$,  $K=0 \ldots \caK$ and 
remove super-indices: for $K$ corresponding a given $\bfi$, we set  $ \hat{B}_K (\bfz) = \hat{B}^{p-2}_{\bfi}(\bfz)$,}   $ {B}_K = \hat{B}_K \circ \varphi_{0,\Gamma_C}^{-1}$  and:
\begin{eqnarray}\label{def:space_mult}
\dis \Lambda^h := Span \{ {B_K(x)}, \quad K=0 \ldots \caK \} . %=\hat{B}^{p-2}_{K}\circ \varphi_{0,\Gamma_C}^{-1}(\bfz)
\end{eqnarray}
For further use, for $v \in \dis L^2(\GC)$ {and for each $K=0 \ldots \caK$}, we denote by $\dis (\Pi_\lambda^h \cdot )_K$ {the following weighted average of $v$}: \begin{eqnarray}\label{def:l2proj_K}
\dis  (\Pi_\lambda^h v)_K = \frac{\dis \intGC v B_K \dG}{\dis \intGC B_K \dG},
\end{eqnarray}
and by $\Pi_\lambda^h$ the global {operator} such as: \begin{eqnarray}\label{def:l2proj}
\dis \Pi_\lambda^h v = \sum_{K = 0}^\caK (\Pi_\lambda^h  v)_K B_K.
\end{eqnarray}
We denote by $L^h$ {the subset} of $W^h$ on which the non-negativity holds only at the control points:$$\dis L^h = \{ \varphi^h \in W^h , \quad (\Pi_\lambda^h  \varphi^h)_K \geq 0\quad \forall K\} .$$ 
%  where $\dis (\Pi_\lambda^h  \varphi^h)_K = \intGC \varphi^h B_K \dG / \intGC B_K \dG$ with $\dis B_K$ basis functions in \eqref{}. We denote by $\Pi_\lambda^h $, the $L^2$-projection on $\Lambda^h$ as follows:  
%\begin{eqnarray}\label{def:l2proj}
%\begin{array}{lccl}
%\dis\Pi_\lambda^h : & \dis L^2(\GC) &\rightarrow &\dis\Lambda^h\\
%& \dis v &\mapsto &\dis \sum_{K = 0}^\caK (\Pi_\lambda^h  v)_K B_K
%\end{array}.
%\end{eqnarray}

\noindent We note that $L^h$ is a convex subset of $W^h$. 

Next, we define the discrete space of the Lagrange multipliers as the negative cones of $L^h$ by $$\dis M^h :=L^{h,*} = \{ \mu^h \in \Lambda^h, \quad \intGC \mu^h \varphi^h \dG \leq 0 \quad \forall \varphi^h \in L^{h} \}.$$

\noindent For any $Q_C \in \restr{\caQ_h}{\GC}$, $\tilde{Q}_C$ denotes the support extension of $Q_C$ (see \cite{daveiga06,daveiga2014}) defined as the image of supports of B-Splines that are not zero on $\hat{Q}_C= \varphi_{0,\Gamma_C}^{-1} (Q_C)$.

\noindent We notice that the operator verifies the following estimate error:
\begin{lemma}\label{lem:op_h1}
Let $\psi \in H^s(\Gamma_C)$ with $0\leq s\leq 1$, the estimate for the local interpolation error reads: 
 \begin{eqnarray} \label{ineq:local disp}
 \dis  \norm{\psi  -\Pi_\lambda^h (\psi )}_{0,\GCQ} \lesssim h^s \norm{\psi}_{s,\GCQtilde} , \qquad \forall \GCQ \in \restr{\caQ_h}{\GC}.
\end{eqnarray}
%where $\GCQtilde$ demotes the support extension of $\GCQ$ relative to the degree $p-2$.
\end{lemma}
\noindent \pf First, Let $c$ be a constant. It holds:
\begin{eqnarray}
\begin{array}{lcl}
\dis \Pi_\lambda^h c &=& \dis \sum_{K = 0}^\caK (\Pi_\lambda^h  c)_K B_K =  \sum_{K = 0}^\caK \frac{\intGC c B_K \dG}{ \intGC B_K \dG}  B_K  =  \sum_{K = 0}^\caK c \frac{\intGC B_K \dG}{ \intGC B_K \dG}  B_K = c \sum_{K = 0}^\caK B_K .
\end{array}\nonumber
\end{eqnarray}
{Using that B-Splines are a partition of the unity, we obtain $ \Pi_\lambda^h c =c$.}

\noindent Let $\psi \in H^s(\Gamma_C)$, it holds:
\begin{eqnarray}\label{ineq:v1}
\begin{array}{lcl}
  \dis \norm{\psi  -\Pi_\lambda^h (\psi )}_{0,\GCQ} &=& \dis \norm{\psi-c  -\Pi_\lambda^h (\psi -c)}_{0,\GCQ}\\[0.2cm]
  & \leq & \dis \norm{\psi-c}_{0,\GCQ} + \norm{\Pi_\lambda^h (\psi -c)}_{0,\GCQ}\\[0.2cm]
 % & \leq & \dis \norm{\psi-c}_{0,\GCQ} + \norm{\Pi_\lambda^h} _{\caL\left(L^2(\GCQtilde);L^2(\GCQ)\right)}\norm{\psi -c}_{0,\GCQtilde}.
\end{array}
\end{eqnarray}
We need now to bound the operator $\Pi_\lambda^h$. We obtain:
\begin{eqnarray}
\begin{array}{lcl}
\dis \norm{\Pi_\lambda^h (\psi -c )}_{0,\GCQ} &=& \dis \norm{ \sum_{K = 0}^\caK \frac{\intGC (\psi -c ) B_K \dG}{ \intGC B_K \dG}  B_K}_{0,\GCQ} \\[0.4cm]
&\leq& \dis \sum_{K: \ supp B_K \cap Q_C \neq \emptyset}^\caK \abs{  \frac{\intGC (\psi -c ) B_K \dG}{ \intGC B_K \dG} }  \norm{B_K}_{0,\GCQ} \\[0.4cm]
&\leq& \dis \sum_{K: \ supp B_K \cap Q_C \neq \emptyset}^\caK  \norm{(\psi -c )}_{0,\GCQtilde}  \frac{\norm{B_K}_{0,\GCQtilde}}{\intGC B_K \dG} \norm{B_K}_{0,\GCQ}.
\end{array}\nonumber
\end{eqnarray}
Using $\norm{B_K}_{0,\GCQtilde} \sim \abs{ \GCQtilde}^{1/2}$, $\norm{B_K}_{0,\GCQ} \sim \abs{{\GCQ}}^{1/2}$, $\dis \intGC B_K \dG \sim \abs{\GCQtilde}$ and Assumption 1, it holds:
\begin{eqnarray}\label{ineq:v2}
\begin{array}{lcl}
\dis \norm{\Pi_\lambda^h (\psi -c)}_{0,\GCQ} &\lesssim& \dis  \norm{(\psi -c )}_{0,\GCQtilde}.
\end{array}
\end{eqnarray}
Using the previous inequalities \eqref{ineq:v1} and \eqref{ineq:v2}, for $0 \leq s \leq 1$, we obtain:
\begin{eqnarray}
\begin{array}{lcl}
  \dis \norm{\psi  -\Pi_\lambda^h (\psi )}_{0,\GCQ}& \lesssim & \dis  \norm{\psi -c}_{0,\GCQtilde}%\\
 % &\lesssim&
\ \lesssim \dis h^s_{\GCQtilde}  \abs{\psi}_{s,\GCQtilde}
\end{array}\nonumber
\end{eqnarray}

\cqfd \\
{\begin{proposition}\label{prop:infsupL2}
For $h$ sufficiently small, there exists a $\beta>0$ such that: 
\begin{eqnarray} \label{ineq:infsupL2}
\dis  \inf_{\mu^h \in M^h}  \sup_{\psi^h  \in  W^h} \frac{ - \intGC  \psi^h \mu^h \dG}{\norm{\psi^h}_{0,\GC}\norm{\mu^h}_{0,\GC}  } \geq \beta .
\end{eqnarray} 
\end{proposition} }
 \noindent \pf In the article \cite{brivadis15}, the authors prove that, if $h$ is sufficiently small, there exists a constant $\beta$ independent of $h$ such that: 
\begin{eqnarray} \label{ineq:infsup_ericka}
\dis \forall {\phi^h} \in (\Lambda^h)^d , \quad \exists u^h \in \restr{V^h}{\GC}, \quad \textrm{s.t.} \quad \dis \frac{ - \intGC {\phi^h} \cdot  u^h \dG}{\norm{u^h}_{0,\Gamma_C} } \geq \beta \norm{{\phi^h}}_{0,\GC} .
\end{eqnarray} 

 \noindent Given now a $\lambda^h \in \Lambda^h$ and {$\psi^h \in W^h$}, we should like to choose  ${\phi^h} =  {\lambda^h n}$ and {$\psi^h = u^h \cdot n$} in \eqref{ineq:infsup_ericka}, but, unfortunately, it is clear that  ${\phi^h} \not\in (\Lambda^h)^d$. Indeed, \eqref{ineq:infsupL2} can obtained from \eqref{ineq:infsup_ericka} via a superconvergence argument that we discuss in the next lines. 

\noindent Let $\Pi_{(\Lambda^h)^d} : L^2(\GC)^d \rightarrow (\Lambda^h)^d$ be a quasi-interpolant defined and studied in \textit{e.g.} see \cite{daveiga2014}.
%%%%%%%% (see \cite{daveiga06,daveiga2014} for more details).
%%%%%%%% (the Schumaker's quasi-interpolant).
%It implies:
%\begin{eqnarray}
%\begin{array}{lcl}
%\dis \norm{\ulambda^h \cdot n - \Pi_{(\Lambda^h)^d}(\ulambda^h \cdot n)}_{0,\Gamma_C}&\leq&\dis C h^{p-1} \abs{\ulambda^h \cdot n }_{p-1,\Gamma_C} \\
%&\leq&\dis C h^{p-1} \sum_{s=0}^{p-1} \norm{ \nabla^s\ulambda^h \cdot \nabla^{p-1-s} n }_{0,\GC}.
%\end{array}\nonumber
%\end{eqnarray}
%If $s=p-1$, then we get $\nabla^s\ulambda^h =0$. 

\noindent If $n \in W^{p-1,\infty} (\GC)$, by the same super-convergence argument used in \cite{brivadis15}, we obtain that:
\begin{eqnarray}\label{ineq:alpha}
\begin{array}{lcl}
\dis \norm{{\phi^h} - \Pi_{(\Lambda^h)^d}({\phi^h})}_{0,\Gamma_C}
%&\leq&\dis C h^{p-1} \sum_{s=0}^{p-2} \norm{ \nabla^s\ulambda^h \cdot \nabla^{p-1-s} n }_{0,\GC}\\[0.6cm]
%&\leq&\dis C  h^{p-1} \abs{\ulambda^h \cdot n }_{p-2,\Gamma_C}\\[0.2cm]
%&\leq& \dis C  h \norm{\ulambda^h \cdot n }_{0,\Gamma_C}\\[0.2cm]
&\leq& \dis \alpha  h \norm{{\phi^h} }_{0,\Gamma_C}.
\end{array}
\end{eqnarray}
Note that: 
\begin{eqnarray}
\begin{array}{lcl}
\dis b(\lambda^h, u^h)& =&\dis - \intGC \lambda^h (u^h \cdot n) \dG=- \intGC {\phi^h} \cdot  u^h \dG \\
& =&\dis - \intGC  \Pi_{(\Lambda^h)^d}({\phi^h}) \cdot  u^h \dG - \intGC  \left({\phi^h} -  \Pi_{(\Lambda^h)^d}({\phi^h}) \right)\cdot  u^h \dG .
\end{array}\nonumber
\end{eqnarray}
By $\inf-\sup$ condition \eqref{ineq:infsup_ericka}, we get:
\begin{eqnarray}
\dis \sup_{u^h \in V^h} \frac{- \intGC  \Pi_{(\Lambda^h)^d}({\phi^h}) \cdot  u^h \dG}{\norm{u^h}_{0,\GC} } \geq \beta \norm{\Pi_{(\Lambda^h)^d}({\phi^h})}_{0,\GC}, \nonumber
\end{eqnarray} 
By \eqref{ineq:alpha}, it holds:
$$\intGC  \left({\phi^h} -  \Pi_{(\Lambda^h)^d}({\phi^h}) \right)\cdot  u^h \dG \leq \alpha h  \norm{{\phi^h} }_{0,\Gamma_C}  \norm{u^h }_{0,\Gamma_C}  .$$
Thus:
\begin{eqnarray}
\begin{array}{lcl}
\dis  \frac{b( \lambda^h ,  u^h )}{~\norm{u^h}_{0,\Gamma_C} } &\geq&\dis  \beta \norm{\Pi_{(\Lambda^h)^d}({\phi^h})}_{0,\GC} - \alpha h  \norm{{\phi^h} }_{0,\Gamma_C} .
\end{array}\nonumber
\end{eqnarray}
Noting that $  \norm{\Pi_{(\Lambda^h)^d}({\phi^h})}_{0,\GC} \geq \norm{{\phi^h} }_{0,\Gamma_C} - \alpha h  \norm{{\phi^h} }_{0,\Gamma_C}$, $ \norm{{\phi^h} }_{0,\Gamma_C} \sim  \norm{\lambda^h }_{0,\Gamma_C}$ and {$ \norm{u^h }_{0,\Gamma_C} \sim  \norm{\psi^h }_{0,\Gamma_C}$}.
Finally, we obtain:
 \begin{eqnarray}
\begin{array}{lcl}
% \dis \sup_{u^h \in V^h} \frac{ b( \lambda^h ,  u^h )}{~\norm{u^h}_{0,\Gamma_C} }&\geq&\dis  \beta \norm{\lambda^h}_{0,\GC} - \alpha h \norm{\lambda^h}_{0,\GC}.
\dis  \sup_{u^h \in V^h}  \frac{ - \intGC  \psi^h \lambda^h \dG}{\norm{\psi^h}_{0,\GC} } &\geq&\dis  \beta \norm{\lambda^h}_{0,\GC} - \alpha h \norm{\lambda^h}_{0,\GC}.
\end{array}\nonumber
\end{eqnarray}

{For $h$ is sufficiently small,} this implies that there exists a constant $\beta'$ independent of $h$ such that:
 \begin{eqnarray}\label{ineq:infsup l2}
\begin{array}{lcl}
% \dis \sup_{u^h \in V^h} \frac{ b( \lambda^h ,  u^h )}{~\norm{u^h}_{0,\Gamma_C} }&\geq&\dis  \beta' \norm{\lambda^h}_{0,\GC}.
 \dis \sup_{u^h \in V^h} \frac{  - \intGC  \psi^h \lambda^h \dG }{~\norm{\psi^h}_{0,\Gamma_C} }&\geq&\dis  \beta' \norm{\lambda^h}_{0,\GC}.
\end{array}
\end{eqnarray} 
\cqfd

{\begin{lemma}\label{lem:neg}
%{Let $\dis \mu^h \in \dis \Lambda^h$ and $\dis \mu^h = \dis \sum_K \mu^h_K B_K$, if we have $\dis \mu^h \in \dis M^h$ then it holds for all $\dis K = 0 \ldots \caK,$ $\mu^h_K \leq 0$.} 
%{It implies that $M^h$ is the negative span of $B_K$.}
For $h$ is sufficiently small, $M^h$ can be characterised as follow: 
 $$M^h \equiv \{  \mu^h \in \Lambda^h, \quad \mu^h = \sum_K \mu^h_K B_K, \quad \mu^h_K\leq 0 \} .$$ 
 \end{lemma} 
 \noindent \pf By definition of $\dis M^h$, we have $\dis \intGC  \mu^h \varphi^h \dG \leq 0$ for all $\dis \varphi^h \in L^h$. We have :
 $$ \dis \intGC  \mu^h \varphi^h \dG =  \sum_K \mu^h_K \intGC B_K \varphi^h \dG = \sum_K \mu^h_K (\Pi_\lambda^h \varphi^h)_K \intGC B_K\dG \leq 0 .$$
 Thus the question becomes if it is possible to find for each $K= 0 \ldots \caK$ a function $\varphi^h_K$ such that:
 $$ \dis \intGC  \mu^h \varphi^h_K \dG = \mu^h_K \intGC B_K\dG. $$
 The existence of such a $\varphi^h_K$ is guaranteed by Proposition \ref{prop:infsupL2}, as a consequence of the $\inf-\sup$ condition \eqref{ineq:infsupL2} is that the matrix representing the scalar product $\dis \intGC  \mu^h \varphi^h \dG$ is full rank. }
  \cqfd \\
 
Then a discretized mixed formulation of the problem \eqref{eq:mixed_form} consists in finding $(u^h,\lambda^h) \in V^h \times M^h$ such that:
\begin{eqnarray} \label{eq:discrete_mixed_form}
\left\{
\begin{array}{rl}
\dis a(u^h,v^h) - b(\lambda^h,v^h) = L(v^h),& \dis \qquad \forall v^h \in V^h,\\
\dis b(\mu^h-\lambda^h,u^h) \geq 0,&\dis \qquad \forall \mu^h \in M^h.
\end{array}
\right.
\end{eqnarray} 
According to Lemma \ref{lem:neg}, we get: $$\{ \mu^h \in M^h: \quad b(\mu^h,v^h)=0 \quad \forall v^h \in V^h) \} =\{0\} ,$$
and using the ellipticity of the bilinear form $a(\cdot,\cdot)$ on $V^h$, then the problem \eqref{eq:discrete_mixed_form} admits a unique solution $(u^h,\lambda^h) \in V^h \times M^h$.

Before addressing the analysis of \eqref{eq:discrete_mixed_form}, let us recall that {the following inequalities} (see \cite{daveiga06}) are true for the primal and the dual space.
\begin{theorem}
Let a given quasi-uniform mesh and let $r, s$ be such that $0 \leq r \leq s \leq p +1$. Then, there exists a constant depending only on $p, \theta, \varphi_0$ and $\hat{W}$ such that for any $v \in H^s(\Omega)$ there exists an approximation $v^h \in V^h$ such that
\begin{eqnarray} \label{ineq:Np}
\dis \norm{v-v^h}_{r,\Omega} \lesssim h^{s-r} \norm{v}_{s,\Omega}.
\end{eqnarray} 
\end{theorem} 
\noindent We will also make use of the local approximation estimates for splines quasi-interpolants that can be found \textit{e.g.} in \cite{daveiga06,daveiga2014}. %We can prove in a similarly way the following lemma. 
\begin{lemma}\label{lem:op_mult}
Let $\lambda \in H^s(\GC)$ with $0\leq s \leq p-1$, then there exists a constant depending only on $p,\varphi_0$ and $\theta$, there exists an approximation $\lambda^h \in \Lambda^h$ such that:
\begin{eqnarray} \label{ineq:local}
\dis h^{-1/2} \norm{\lambda-\lambda^h}_{-1/2,\GCQ} +  \norm{\lambda-\lambda^h}_{0,\GCQ} \lesssim h^{s} \norm{\lambda}_{s,\GCQtilde} , \qquad \forall \GCQ \in \restr{\caQ_h}{\GC}.
\end{eqnarray} 
 %where $h_Q$ is the size of $Q \in \mathcal{Q}_h$.\dis \norm{\lambda-\lambda^h}_{0,\GCQ} \leq C h_Q^{s} \norm{\lambda}_{s,\GCQtilde}\leq C h^{s} \norm{\lambda}_{s,\GCQtilde},
\end{lemma}

{\noindent It is well known \cite{fortinbook13} that the stability for the mixed problem \eqref{eq:mixed_form} is linked to the $\inf-\sup$ condition.
\begin{theorem}
For $h$ sufficiently small, $n$ sufficiently regular and for any $\mu^h  \in \Lambda^h$, it holds:
\begin{eqnarray} \label{ineq:infsup}
\dis \sup_{v^h \in V^h} \frac{b(\mu^h,v^h)}{\norm{v^h}_{V} } \geq \beta \norm{\mu^h}_{W'},
\end{eqnarray} 
where $\beta$ is independent of $h$.
\end{theorem} }

\noindent \pf
{By Lemma \ref{lem:neg},} there exists %a $\Pi$ 
a Fortin's operator $\Pi : \ L^2(\Gamma_C) \rightarrow \restr{V^h}{\GC} \cap H^1_0(\Gamma_C) $ such that 
\begin{eqnarray}
\begin{array}{ll}
&\dis b(\lambda,\Pi (u))= b(\lambda,u) , \ \forall \lambda \in M \quad \textrm{ and } \quad \dis \norm{\Pi (u)}_{0,\Gamma_C} \leq \norm{u}_{0,\Gamma_C}.
\end{array}\nonumber
\end{eqnarray}
Let $I_h$ be a $L^2$ and $H^1$ stable quasi-interpolant onto $\restr{V^h}{\GC}$ (for example, the Schumaker's quasi-interpolant, see for more details \cite{daveiga2014}). It is important to notice that $I^h$ preserves the homogeneous Dirichlet boundary condition.

\noindent We set $\Pi_F = \Pi (I-I_h) + I_h$. It is classical to see that:
\begin{eqnarray}\label{prop:pif-b}
\begin{array}{ll}
 \dis b(\lambda,\Pi_F (u)) =  b(\lambda,u) , \quad  \forall \lambda \in M,
 \end{array}
\end{eqnarray}
%Indeed, by definition of $\Pi$, we obtain: $$\dis b(\lambda,\Pi_F (u)-u) = b(\lambda, \Pi (u -I_hu )+ I_hu -u ) = b(\lambda, \Pi (u -I_hu ))  -b(\lambda, u - I_hu))=0 .$$
and it is easy to see that: 
%At last, it holds: 
\begin{eqnarray}\label{prop:pif-=}
\begin{array}{ll}
\dis \Pi_F (u^h)= u^h , \quad \forall u^h \in \restr{V^h}{\GC} .
\end{array}
\end{eqnarray}
Moreover, by stability of $\Pi$ and $I_h$, it holds: 
\begin{eqnarray}\label{prop:pif-leq}
\begin{array}{ll}
\dis \dis \norm{\Pi_F (u)}_{0,\Gamma_C} \lesssim \norm{u}_{0,\Gamma_C},  \quad \forall u \in L^2(\GC).
\end{array}
\end{eqnarray}
and also
\begin{eqnarray}\label{prop:pif-leqh1}
 \begin{array}{lcl}
 \norm{\Pi_F (u)}_{1,\Gamma_C} \lesssim \norm{ (u)}_{1,\Gamma_C} , \quad \forall u \in H^1(\GC).
  \end{array}
\end{eqnarray}
%Indeed, for $u^h \in \restr{V^h}{\GC}$, since $I_hu^h = u^h$, it implies: $$\dis \Pi_F (u^h) = \Pi (u^h -I_hu^h )+ I_hu^h  = \Pi (u^h -u^h )+ u^h= u^h.$$
%Now, let us prove that:  \begin{eqnarray}\label{prop:pif-leqh1}
% \begin{array}{lcl}
% \norm{\Pi_F (u)}_{1,\Gamma_C} \lesssim \norm{ (u)}_{1,\Gamma_C} , \quad \forall u \in H^1(\GC).
%  \end{array}
%\end{eqnarray}
%Using the discrete norm inequality for a quasi-uniform mesh, the $L^2$-stability of $\Pi$ and the $H^1$-stability of the operator $I_h$, for $u \in H^1(\GC)$, it holds:
% \begin{eqnarray}
% \begin{array}{lcl}
%\dis \dis \norm{\Pi_F (u)}_{1,\Gamma_C} &\leq&\dis \norm{\Pi (u -I_hu )}_{1,\Gamma_C}+ \norm{I_hu}_{1,\Gamma_C}\\[0.2cm]
%&\lesssim&\dis h^{-1} \norm{\Pi (u -I_hu )}_{0,\Gamma_C}+ \norm{I_hu}_{1,\Gamma_C}\\[0.2cm]
%&\lesssim&\dis h^{-1} \norm{ u -I_hu }_{0,\Gamma_C}+ \norm{u}_{1,\Gamma_C}\\[0.2cm]
%&\lesssim&\dis \norm{u}_{1,\Gamma_C}.
%\end{array}\nonumber
%\end{eqnarray}
{To conclude, we distinguish between two cases :
\begin{itemize}
\item If $\overline{\Gamma}_D \cap \overline{\Gamma}_C =\emptyset$, it is well know that $W=H^{1/2}(\GC)$. By interpolation of Sobolev Spaces, using \eqref{prop:pif-leq} and \eqref{prop:pif-leqh1}, we obtain:
$$ \dis  b(\lambda,\Pi_F (u)) = b(\lambda,u), \quad \forall \lambda \in M \quad \textrm{ and } \quad \dis \norm{\Pi_F (u)}_{W} \lesssim \norm{u}_{W}. $$
Then $\inf-\sup$ condition \eqref{ineq:infsup} holds thanks to Proposition 5.4.2 of \cite{fortinbook13}.
\item If $\overline{\Gamma}_D \cap \overline{\Gamma}_C \neq \emptyset$, it is enough to remind that for all $u \in H^1_{0, \Gamma_D \cap \Gamma_C} (\GC)$, we have $\dis \Pi_F (u) \in H^1_{0, \Gamma_D \cap \Gamma_C} (\GC)$ and \eqref{prop:pif-leqh1} is valid on the subspace $H^1_{0, \Gamma_D \cap \Gamma_C} (\GC)$.
Again by interpolation argument between \eqref{prop:pif-leq} and \eqref{prop:pif-leqh1}, it holds $ \dis \norm{\Pi_F (u)}_{W} \leq C\norm{u}_{W}$ which ends the proof.
\end{itemize}}
\cqfd

%%%%%%%%%%%%%%%%%%%%%%%%%%%%%%%%%%%%%%%%%%%%%%%%

%%%%%%%%%%%%%%%%%%%%%%%%%%%%%%%%%%%%%%%%%%%%%%%% 
  
\section{\textit{A priori} error analysis}
\label{sec:sec2}
%\subsection{Unilateral contact problem}
%\label{subsec:continuous_pb}

In this section, we present an {optimal} \textit{a priori} error estimate for the Signorini mixed problem. Our estimates follows the ones for finite elements, provided in \cite{coorevits-hild-lhalouani-sassi-02,hild-laborde-02}, and refined in \cite{drouet15}. In particular, in \cite{drouet15} the authors overcome a technical assumption on the geometric structure of the contact set and we are able to avoid such as assumptions also in our case.
%{Due to the intrinsic lack of regularity of contact solutions, in general, only $p=2$ and $p=3$ are considered. We focus on case $p=2$, for higher order a specific local treatment would be necessary to recover origner rates. }
%The following work is the prolongation to the isogeometric analysis of the classical result \cite{hild-laborde-02}, where an estimate for mixed variational formulation with quadratic finite elements method. This estimate is achieved by the use of an assumption of the set contact. A recent work \cite{drouet15} passes over this assumption in two and three dimensions with linear and quadratic elements thanks to the definition of a contact set and a non-contact set.
%{Due to the intrinsic lack of regularity of contact solutions, we limit our study to the case $p=2$ and $p=3$ and for this cases, and in the small deformation regime, we are able to provide a full convergence analysis for the strategy we propose. Indeed, following \cite{drouet15,coorevits-hild-lhalouani-sassi-02}, we are able to prove that our method shows optimal convergence properties (despite the lower degree in multipliers) for general contact solution having Sobolev regularity up to $5/2$ (\cite{Moussaoui92}). }

Indeed, for any $p$, we prove our method to be optimal for solutions with regularity up to $5/2$. 
{Thus, optimality for the displacement is obtained for any $p \geq 2$.} The cheapest {and more convenient} method proved optimal corresponds to the choice $p=2$. Larger values of $p$ may be of interest because they produce continuous pressures, but, on the other hand, the error bounds remain limited by the regularity of the solution, \textit{i.e.} , up to $C h^{3/2}$. Clearly, to enhance approximation suitable local refinement may be used, \cite{lorenzis11,lorenzis12}, but this choice outside the scope of this paper.

In order to prove Theorem \ref{tmp:apriori} which follows, we need a few preparatory Lemmas.

\noindent First, we introduce some notation and some basic estimates. Let us define the active-set strategy for the variational problem. Given an element $Q_C \in \restr{\caQ_h}{\GC}$ of the undeformed mesh, we denote by $Z_{C}(\GCQ)$ the contact set and by $Z_{NC}(\GCQ)$ the non-contact set in $\GCQ$, as follows: 
\begin{eqnarray} %\label{}
\begin{array}{rl}
\dis  Z_{C}(\GCQ) = \{ x \in \GCQ, \quad u_n (x) = 0\}  \quad \textrm{and} \quad \dis  Z_{NC}(\GCQ) = \{ x \in \GCQ, \quad u_n (x) > 0\} .
\end{array}\nonumber
\end{eqnarray}
$|Z_{C}(\GCQ)|$ and $|Z_{NC}(\GCQ)|$ stand for their measures and $|Z_{C}(\GCQ)|+|Z_{NC}(\GCQ)| = |Q_C| = C h_{Q_C}^{d-1}$.
\begin{remark}\label{rq:cont}
Since $u_n$ belongs to $H^{1+\nu}(\Omega)^{2}$ for $0<\nu<1$, if $d=2$ the Sobolev embeddings ensure that $u_n \in \caC^0(\partial \Omega)$. It implies that $Z_{C}(\GCQ)$ and  $Z_{NC}(\GCQ)$ are measurable as inverse images of a set by a continuous function.
\end{remark}
The following estimates are the generalization to the mixed problem of Lemma 2 of Appendix of the article \cite{drouet15}. We recall that if $(u,\lambda)$ is a solution of the mixed problem \eqref{eq:mixed_form} then $\sigma_n(u)=\lambda$. So, the following lemma can be proven exactly in the same way.
\begin{lemma} \label{lemma1}
Let $d=2$ or $3$. Let $(u,\lambda)$ be the solution of the mixed formulation \eqref{eq:mixed_form} and let $u \in H^{3/2+ \nu}(\Omega)^d$ with $0<\nu<1$. Let $h_Q$ the be the diameter of the trace element $\GCQ$ and the set of contact $Z_{C}(\GCQ)$ and non-contact $Z_{NC}(\GCQ)$ defined previously in $\GCQ$. \\
We assume that $|Z_{NC}(\GCQ)|>0$, the following %$L^1$- and 
$L^2$-estimate holds for $\lambda$:
%\begin{eqnarray} \label{estimate1}
%\dis \norm{\lambda}_{L^1(\GCQ)} \leq \frac{|Z_{C}|^{1/2}}{|Z_{NC}|^{1/2}} h_Q^{d/2+\nu-1/2} |\lambda|_{\nu,\GCQ } ,
%\end{eqnarray}\\[-0.9cm]
\begin{eqnarray} \label{estimate2}
\dis \norm{\lambda}_{0,\GCQ} \leq \frac{1}{|Z_{NC}(\GCQ)|^{1/2}} h_{Q_C}^{d/2+\nu-1/2} |\lambda|_{\nu,\GCQ } .
\end{eqnarray}
We assume that $|Z_{C}(\GCQ)|>0$, the following %$L^1$- and 
$L^2$-estimates hold for $\nabla u_n$:
%\begin{eqnarray} \label{estimate3}
%\dis \norm{\nabla u_n}_{L^1(\GCQ)} \leq \frac{|Z_{NC}|^{1/2}}{|Z_{C}|^{1/2}} h_Q^{d/2+\nu-1/2} |\nabla u_n|_{\nu,\GCQ } ,
%\end{eqnarray}\\[-0.9cm]
\begin{eqnarray} \label{estimate4}
\dis \norm{\nabla u_n}_{0,\GCQ} \leq \frac{1}{|Z_{C}(\GCQ)|^{1/2}} h_{Q_C}^{d/2+\nu-1/2} |\nabla u_n|_{\nu,\GCQ } .
\end{eqnarray}
\end{lemma}

\begin{theorem}\label{tmp:apriori} Let $(u,\lambda)$ and $(u^h,\lambda^h)$ be respectively the solution of the mixed problem \eqref{eq:mixed_form} and the discrete mixed problem \eqref{eq:discrete_mixed_form}. Assume that $u \in H^{3/2+\nu}(\Omega)^d$ with $0<\nu<1$. Then, the following error estimate is satisfied:
\begin{eqnarray} \label{eq:apriori}
\norm{u-u^h}_{V}^2 + \norm{\lambda-\lambda^h}^2_{W'} \lesssim h^{1+2\nu} \norm{u}^2_{3/2 + \nu,\Omega}.
\end{eqnarray} 

%\begin{remark}
%In two-dimensional case
%\end{remark}

\end{theorem}
\noindent \pf {In the article \cite{hild-laborde-02} Proposition 4.1, it is proved that if $(u,\lambda)$ is the solution of the mixed problem \eqref{eq:mixed_form} and $(u^h,\lambda^h)$ is the solution of the discrete mixed problem \eqref{eq:discrete_mixed_form}, it holds :} 

\begin{eqnarray} %\label{}
\begin{array}{rl}
\dis \norm{u-u^h}_{V}^2 +  \norm{\lambda-\lambda^h}_{W'}^2 \lesssim & \dis %\inf_{v^h \in N^{p}}
\norm{u-v^h}_{V}^2 + %\inf_{\mu^h \in \Lambda^h}
\norm{\lambda-\mu^h}_{W'}^2 \\[0.4cm]
& \dis +\max(-b(\lambda,u^h),0)  +\max(- b(\lambda^h,u),0) .
\end{array}\nonumber
\end{eqnarray}
It remains to estimate on the previous inequality the last two terms to obtain the estimate \eqref{eq:apriori}. \\

\noindent\textbf{Step 1: estimate of $-b(\lambda,u^h) \dis =   \intGC \lambda u^h_n  \dG$.}\\
\noindent Using the operator $\Pi_\lambda^h$ defined in \eqref{def:l2proj}, it holds:
\begin{eqnarray} %\label{}
\begin{array}{lcl}
\dis -b(\lambda,u^h) & \dis =&\dis \intGC \lambda u^h_n \dG \dis= \intGC \lambda  \left(u^h_n - \Pi_\lambda^h (u^h_n)\right) \dG + \intGC \lambda \Pi_\lambda^h ( u^h_n) \dG \\[0.4cm]
& \dis = &\dis  \intGC \left(\lambda -\Pi_\lambda^h (\lambda) \right) \left(u^h_n - \Pi_\lambda^h (u^h_n)\right) \dG + \intGC \Pi_\lambda^h (\lambda) \left(u^h_n - \Pi_\lambda^h (u^h_n)\right) \dG \\[0.4cm]
&& \dis + \intGC \lambda \Pi_\lambda^h ( u^h_n) \dG. 
\end{array}\nonumber
\end{eqnarray}
\noindent Since $\lambda$ is a solution of \eqref{eq:mixed_form}, it holds $\dis \Pi_\lambda^h (\lambda) \leq 0$. Furthermore, $u^h$ is a solution of \eqref{eq:discrete_mixed_form}, thus $\dis \intGC \Pi_\lambda^h (\lambda) \left(u^h_n - \Pi_\lambda^h (u^h_n)\right) \dG \leq 0$ and $\dis \intGC \lambda \Pi_\lambda^h ( u^h_n) \dG \leq 0$.

\noindent We obtain:
\begin{eqnarray} \label{eq:eq1}
\begin{array}{lcl}
\dis -b(\lambda,u^h)& \dis \leq& \dis   \intGC \left(\lambda -\Pi_\lambda^h (\lambda) \right) \left(u^h_n - \Pi_\lambda^h (u^h_n)\right) \dG \\[0.4cm]
&  \leq &\dis \intGC \left(\lambda -\Pi_\lambda^h (\lambda) \right) \left(u^h_n -u_n- \Pi_\lambda^h (u^h_n-u_n)\right) \dG \\[0.4cm]
& &+\dis\intGC \left(\lambda -\Pi_\lambda^h (\lambda) \right) \left(u_n- \Pi_\lambda^h (u_n)\right) \dG.
\end{array}
\end{eqnarray}
The first term of \eqref{eq:eq1} is bounded in an optimal way by using \eqref{ineq:local disp}, the summation on each physical element, Theorem \ref{thm:u-lambda} and the trace theorem:
\begin{eqnarray} %\label{}
\begin{array}{ll}
\dis \intGC \left(\lambda -\Pi_\lambda^h (\lambda) \right) \left(u^h_n -u_n- \Pi_\lambda^h (u^h_n-u_n)\right) \dG \!\!\!\!\!&\leq \dis \norm{\lambda -\Pi_\lambda^h (\lambda) }_{0,\GC} \norm{u^h_n -u_n- \Pi_\lambda^h (u^h_n-u_n)}_{0,\GC}\\[0.2cm]
&\leq \dis C h^{1/2+\nu} \norm{\lambda }_{\nu,\GC}  \norm{u_n-u^h_n}_{W}\\[0.2cm]
&\leq \dis C h^{1/2+\nu} \norm{u}_{3/2+\nu,\Omega}  \norm{u-u^h}_{V}.
\end{array}\nonumber
\end{eqnarray}
We need now to bound the second term in \eqref{eq:eq1}. Let $Q_C$ be an element of $ \restr{\caQ_h}{\GC}${. I}f either $|\ZC|$ or $|\ZNC|$ are null, the integral on $\GCQ$ vanishes. So we suppose that either $|\ZC|$ or $|\ZNC|$ are greater than {$\dis|\GCQ|/2= C h_{Q_C}^{d-1}$} and we consider the two cases, separately. \\

{Similarly to the article \cite{hild-laborde-02}, we can prove that if:}
\begin{itemize}
\item $|\ZC| \geq \dis {|\GCQ|}/{2} $. Using the estimate \eqref{ineq:local disp}, the estimate \eqref{estimate4} of Lemma \ref{lemma1} and the Young's inequality, it holds: {
\begin{eqnarray} %\label{}
\begin{array}{rcl}
\dis \intGCQ (\lambda -\Pi_\lambda^h (\lambda) ) (u_n  -\Pi_\lambda^h (u_n ) ) \dG 
%&\leq& {begin \ cut \ } \dis \intGCQ \lambda (u_n  -\Pi_\lambda^h (u_n ) ) \dG \\
% &\leq&\dis  \norm{\lambda}_{0,\GCQ}  \norm{u_n  -\Pi_\lambda^h (u_n )}_{0,\GCQ}\\[0.2cm]
% &\leq&\dis C \norm{\lambda}_{0,\GCQ} h \norm{ u_n}_{1,\GCQtilde}\\[0.2cm]
% &\leq& \dis \frac{C}{|\ZC|^{1/2}} h^{d/2+\nu -1/2} \norm{\lambda}_{\nu,\GCQ}  h^{1+\nu}\norm{ u_n}_{1+\nu,\GCQtilde} \\[0.4cm]
% &\leq& \dis C h^{d/2+2\nu +1/2} h^{ -d/2 +1/2} \norm{\lambda}_{\nu,\GCQ} \norm{ u_n}_{1+\nu,\GCQtilde} \\[0.2cm]
% &\lesssim& \dis  h^{1+2\nu} \norm{\lambda}_{\nu,\GCQ} \norm{ u_n}_{1+\nu,\GCQtilde}\\[0.2cm]
&\lesssim&\dis  h^{1+2\nu} (\norm{\lambda}^2_{\nu,\GCQ} + \norm{ u_n}^2_{1+\nu,\GCQtilde}).
\end{array}\nonumber
\end{eqnarray} }

\item $|\ZNC| \geq \dis {|\GCQ|}/{2} $. Using the estimate \eqref{ineq:local disp}% and \eqref{ineq:local}
, the estimate \eqref{estimate2} of Lemma \ref{lemma1} and the Young's inequality, it holds:
\begin{eqnarray} %\label{}
\begin{array}{rcl}
 { \dis \intGCQ (\lambda -\Pi_\lambda^h (\lambda) ) (u_n  -\Pi_\lambda^h (u_n ) ) \dG }
%&\leq& \dis {begin \ cut \ } \norm{\lambda -\Pi_\lambda^h (\lambda)}_{0,\GCQ}  \norm{u_n  -\Pi_\lambda^h (u_n )}_{0,\GCQ} \\[0.4cm]
%%&\leq& \dis C \norm{\lambda}_{0,\GCQtilde} h^{1} \norm{\nabla u_n}_{0,\GCQtilde} \\[0.4cm]
%&\leq& \dis \frac{C}{|\ZNC|^{1/2}} h^{d/2+2\nu +1/2} \norm{\lambda}_{\nu,\GCQtilde} \norm{ u_n}_{1+\nu,\GCQtilde} \\[0.4cm]
%&\lesssim& \dis  h^{1+2\nu} \norm{\lambda}_{\nu,\GCQtilde} \norm{ u_n}_{1+\nu,\GCQtilde}\\[0.2cm]
&{ \lesssim }&{\dis  h^{1+2\nu} (\norm{\lambda}^2_{\nu,\GCQtilde} + \norm{ u_n}^2_{1+\nu,\GCQtilde}). }
\end{array}\nonumber
\end{eqnarray}
\end{itemize}  

\noindent Summing over all the contact elements and {distinguishing the two cases $\ZC \geq {|\GCQ|}/{2}$ and $\ZNC \geq {|\GCQ|}/{2}$}, it holds:
\begin{eqnarray} %\label{}
\begin{array}{rcl}
\dis \intGC (\lambda -\Pi_\lambda^h (\lambda) ) (u_n  -\Pi_\lambda^h (u_n ) ) \dG & = &\dis \sum_{{Q_C \in  \restr{\caQ_h}{\GC}}} \intGCQ (\lambda -\Pi_\lambda^h (\lambda) ) (u_n  -\Pi_\lambda^h (u_n ) ) \dG \\[0.4cm]
&\leq&\dis C h^{1+2\nu} \sum_{Q_C \in  \restr{\caQ_h}{\GC}} \norm{\lambda}^2_{\nu,\GCQ} + \norm{\lambda}^2_{\nu,\GCQtilde} + \norm{ u_n}^2_{1+\nu,\GCQtilde} \\[0.4cm]
&\leq&\dis C h^{1+2\nu} \sum_{Q_C \in  \restr{\caQ_h}{\GC}} \norm{\lambda}^2_{\nu,\GCQ} + \sum_{Q_C' \in \tilde{Q}_C} \norm{\lambda}^2_{\nu,\GCQ'} + \norm{ u_n}^2_{1+\nu,\GCQ'} \\[0.4cm]
&\leq&\dis C h^{1+2\nu} \Big{(} \norm{\lambda}^2_{\nu,\GC} +  \sum_{Q \in  \restr{\caQ_h}{\GC}} \sum_{Q_C' \in \tilde{Q}_C} \norm{\lambda}^2_{\nu,\GCQ'} + \norm{ u_n}^2_{1+\nu,\GCQ'}  \Big{)}. 
\end{array}\nonumber
\end{eqnarray}
Due to the compact supports of the B-Splines basis functions, there exists a constant $C$ depending only on the degree $p$ and the dimension $d$ of the undeformed domain such that: $$\dis   \sum_{Q \in  \restr{\caQ_h}{\GC}} \sum_{Q_C' \in \tilde{Q}_C} \norm{\lambda}^2_{\nu,\GCQ'} + \norm{ u_n}^2_{1+\nu,\GCQ'}  \leq C \norm{\lambda}^2_{\nu,\GC} + C\norm{ u_n}^2_{1+\nu,\GC} .$$
So we have: $$\dis \intGC (\lambda -\Pi_\lambda^h (\lambda) ) (u_n  -\Pi_\lambda^h (u_n ) ) \dG \leq C h^{1+2\nu } \Big{(}  \norm{\lambda}^2_{\nu,\GC} + \norm{ u_n}^2_{1+\nu,\GC}\Big{)} ,$$
\textit{i.e.} $$\dis \intGC (\lambda -\Pi_\lambda^h (\lambda) ) (u_n  -\Pi_\lambda^h (u_n ) ) \dG \leq C h^{1+2\nu}  \norm{u}^2_{3/2+\nu,\Omega} .$$

\noindent We conclude that:
\begin{eqnarray}
\begin{array}{rl}
\dis -b(\lambda,u^h) & \dis \lesssim  h^{1/2+\nu} \norm{u}_{3/2+\nu,\Omega} \norm{u-u^h}_{V} + h^{1+2\nu}\norm{u}^2_{3/2+\nu,\Omega}.
\end{array}\nonumber
\end{eqnarray}

\noindent Using Young's inequality, we obtain:
\begin{eqnarray} \label{step2}
\begin{array}{rl}
-b(\lambda,u^h) \lesssim h^{1+2\nu}  \norm{u}^2_{3/2+\nu,\Omega} + \norm{u-u^h}^2_{V} .
\end{array}
\end{eqnarray}

\noindent\textbf{Step 2: estimate of $-b(\lambda^h,u) \dis =   \intGC \lambda^h u_n  \dG$.}\\
Let us denote by $j^h$ the Lagrange interpolation operator of order one on {$ \restr{\caQ_h}{\GC}$}. 
\begin{eqnarray} %\label{}
\begin{array}{rl}
\dis -b(\lambda^h,u) & \dis =   \intGC \lambda^h u_n  \dG = \intGC \lambda^h (u_n  -j^h (u_n )) \dG + \intGC \lambda^h j^h (u_n ) \dG .\\
\end{array}\nonumber
\end{eqnarray}
Note that by remark \ref{rq:cont}, $u_n$ is continuous and $j^h (u_n )$ is well define. \\
\noindent Since $u$ is a solution of \eqref{eq:mixed_form}, it holds $\dis j^h (u_n ) \geq 0$. Thus, $\dis \intGC \lambda^h j^h (u_n ) \dG \leq 0, \quad \lambda^h \in M^h$.

\noindent As previously, we obtain:
\begin{eqnarray} %\label{}
\begin{array}{rl}
\dis -b(\lambda^h,u) & \dis \leq \intGC \lambda^h u_n  \dG \dis \leq \intGC \lambda^h (u_n  -j^h (u_n )) \dG \\
& \dis \leq \intGC (\lambda^h - \lambda)(u_n  -j^h (u_n )) \dG + \intGC \lambda (u_n  -j^h (u_n )) \dG\\
& \dis \leq \intGC \lambda (u_n  -j^h (u_n )) \dG+  \norm{\lambda-\lambda^h}_{W'}  \norm{u_n  - j^h (u_n )}_{W}\\
& \dis  \leq \intGC \lambda (u_n  -j^h (u_n )) \dG + h^{1/2 +\nu}\norm{u_n }_{1+\nu,\GC} \norm{\lambda-\lambda^h}_{W'}\\
& \dis  \leq \intGC \lambda (u_n  -j^h (u_n )) \dG + h^{1/2 +\nu}\norm{u}_{3/2+\nu,\Omega} \norm{\lambda-\lambda^h}_{W'}.
\end{array}\nonumber
\end{eqnarray}

\noindent Now, we need to show that: \begin{eqnarray} \label{step1_1}
\intGC \lambda (u_n  -j^h (u_n )) \dG \leq C h^{1+2\nu} \norm{u}_{3/2+\nu,\Omega}^2.\end{eqnarray}
 
\noindent The proof of this inequality is done in the paper \cite{drouet15} for both linear and quadratic finite elements, and can be repeated here verbatim. In this proof, two cases are considered:
\begin{enumerate}
\item either $|\ZC|$ or $|\ZNC|$ is null and thus the inequality is trivial;
\item where either $|\ZC|$ or $|\ZNC|$ is greater than $\dis {|\GCQ|}/{2}= C h_{Q_C}^{d-1}$.
\end{enumerate}
 As previously, choosing either $|\ZC|$ or $|\ZNC|$, using the previous Lemma \ref{lemma1} and by summation on all element of mesh, 
{\noindent we conclude that:
\begin{eqnarray} %\label{step1}
\begin{array}{rl}
\dis -b(\lambda^h,u) & \dis  \leq \intGC \lambda (u_n  -j^h (u_n )) \dG + h^{1/2 +\nu}\norm{u}_{3/2+\nu,\Omega} \norm{\lambda-\lambda^h}_{W'}\\
& \dis  \lesssim h^{1+2\nu}\norm{u}^2_{3/2+\nu,\Omega} + h^{1/2 +\nu}\norm{u}_{3/2+\nu,\Omega} \norm{\lambda-\lambda^h}_{W'}.
\end{array}\nonumber
\end{eqnarray}
Using Young's inequality, we obtain:
\begin{eqnarray} \label{step1}
\begin{array}{rl}
\dis -b(\lambda^h,u) & \dis  \lesssim h^{1+2\nu}\norm{u}^2_{3/2+\nu,\Omega} + \norm{\lambda-\lambda^h}^2_{W'}.
\end{array}
\end{eqnarray}}
Finally, we can conclude using \eqref{step1} and  \eqref{step2}.%, we obtain the \textit{a priori} error estimation  \eqref{eq:apriori}. 

\cqfd

%%%%%%%%%%%%%%%%%%%%%%%%%%%%%%%%%%%%%%%%%%%%%%%%

\section{Numerical Study}
\label{sec:sec3}

{In this section, we perform  {a} numerical validation for the method we propose {d} in small as well as in large deformation {frameworks}, \textit{i.e.}, also beyond the theory developed in previous Sections. Due to the intrinsic lack of regularity of contact solutions, we restrict ourselves to the case $p=2$,  {for which the $N_2 / S_0$ method is tested}.}

%In this section, we consider the case $p=2$. Indeed, It should be extended to the case $p=3$, to obtain a multiplier continuous. This section present an {optimal} \textit{a priori} error estimate and from the theoretical point of view, the choice $p=3$ has no advantages for the unilateral contact problem where the maximal regularity is $5/2$. As we said previously, It is well know (see \cite{Moussaoui92}), for unilateral contact problems, the regularity of the solution can generally not be passed beyond $5/2$. Therefore, we analyse the numerical performance for the $N_2 / S_0$ method:
%quadratic NURBS basis functions are used for discretizing displacement
%unknowns and Lagrange multipliers are interpolated by means of piecewise constant functions.
%Below some examples that validate and prove the accuracy of the proposed methods are presented.
The suite of benchmarks reproduces the classical Hertz contact problem \cite{Hertz1882,johnsonbook}:
Sections \ref{subsec:Hertz problems} and \ref{subsec:Hertz problems} analyse
the two and three-dimensional cases for a small deformation setting, whereas
Section \ref{subsec:large deformation} considers the large deformation problem in 2D.
The examples were performed using an in-house code based on the igatools library (see \cite{Pauletti2015} for further details).

{In the following example, to prevent that the {contact zone} is empty, we considered, only for the initial gap, that {there exists contact} if the $\dis g_n\leq 10^{-9}$.}\\

\subsection{Two-dimensional Hertz problem}
\label{subsec:Hertz problems}
\begin{figure}[!ht]
  \centering
  \subfigure[Stress magnitude distribution for the undeformed mesh.]{\label{fig:mesh1_a} \includegraphics[width= 8.5cm]{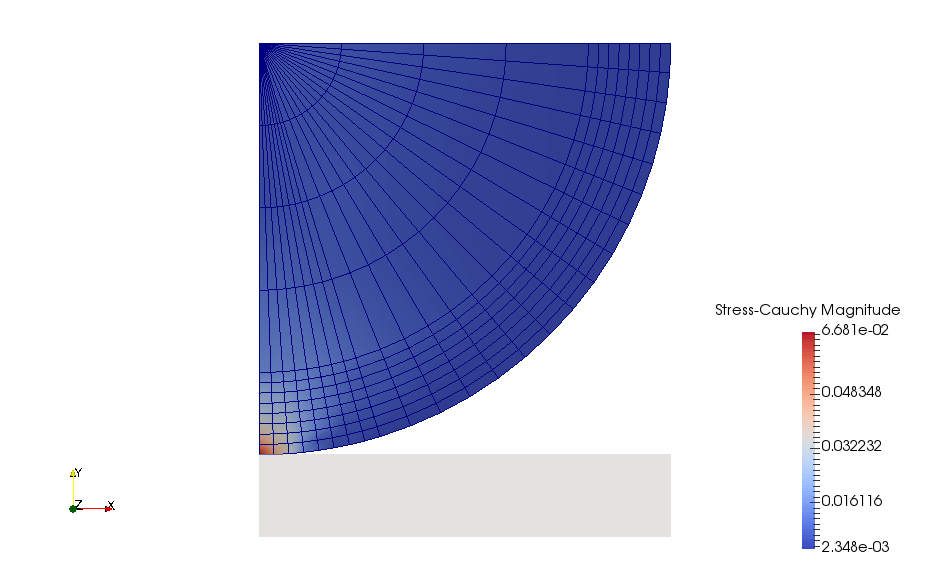}}\hfill
  \subfigure[Analytical and numerical contact pressure.]{\label{fig:mesh1_b} \includegraphics[width= 7cm]{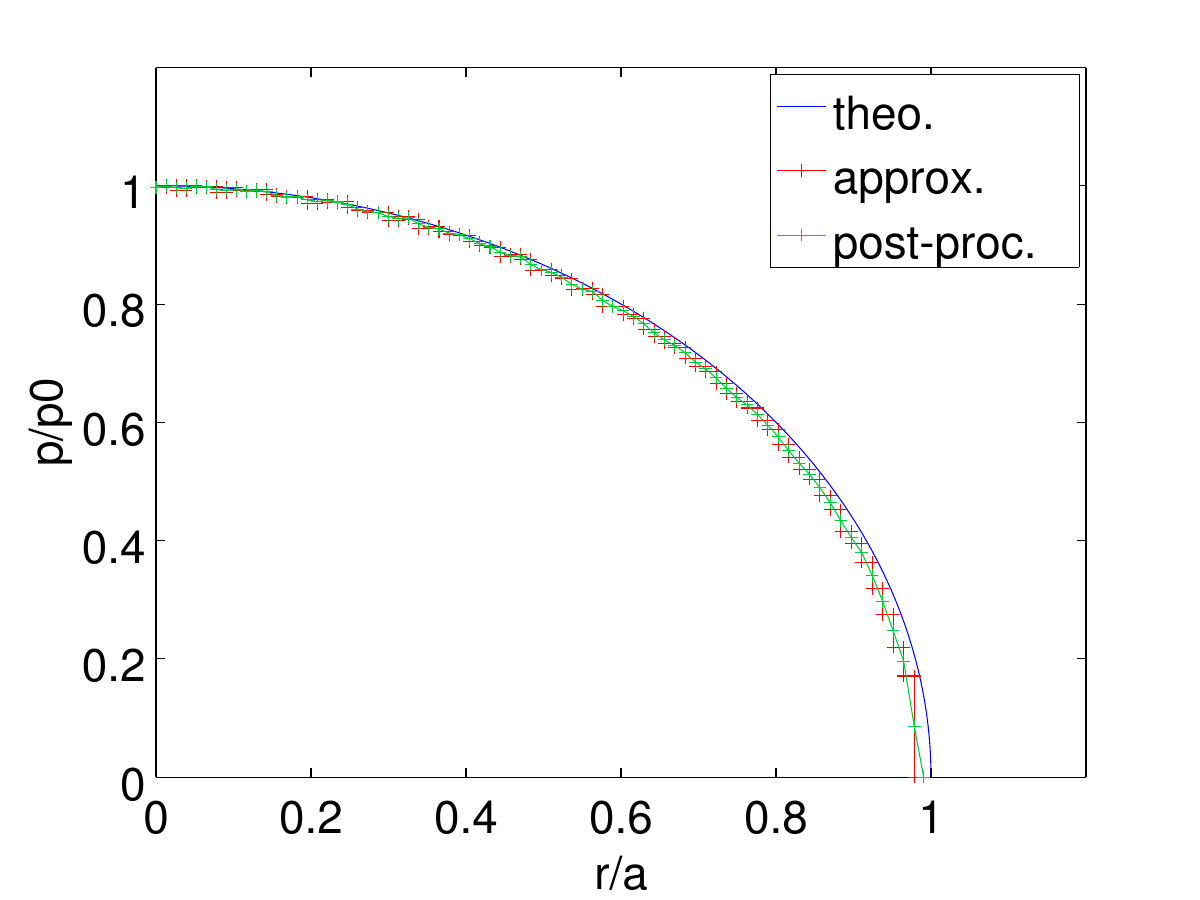}}
  \caption{2D Hertz contact problem with $N_2 / S_0$ method for an applied pressure $P=0.003$.}
  \label{fig:mesh1}
\end{figure}

The first example included in this section analyses the two-dimensional frictionless Hertz contact problem considering small elastic deformations.
It consists in an infinitely long half cylinder body with radius $R=1$, that it is deformable and whose
material is linear elastic, with Young's modulus $E = 1$ and Poisson's ratio $\nu = 0.3$.
A uniform pressure $P=0.003$ is applied on the top face of the cylinder while the curved surface contacts against a horizontal rigid plane
(see Figure \ref{fig:mesh1_a}).
Taking into account the test symmetry and the ideally infinite length of the cylinder,
the problem is modelled as 2D quarter of disc with proper boundary conditions.

Under the hypothesis that the contact area is small compared to the cylinder dimensions, the Hertz's analytical solution (see \cite{Hertz1882,johnsonbook})
predicts that the contact region is an infinitely long band whose width is $2a$, being $a = \sqrt{8R^2P(1- \nu^2)/\pi E}$.
Thus, the normal pressure, {that} follows an elliptical distribution along the width direction $r$, is $p(r)=p_0\sqrt{1-r^2/a^2}$,
where the maximum pressure, at the central line of the band ($r=0$), is $p_0 = 4 RP/ \pi a$.
For the geometrical, material and load data chosen in this numerical test,
the characteristic values of the solution are $a=0.083378$ and $p_0 = 0.045812$.
Notice that, as required by Hertz's theory hypotheses, $a$ is sufficiently small compared to $R$.

It is important to remark that, despite the fact that Hertz's theory provides a
full description of the contact area and the normal contact pressure
in the region, it does not describe analytically the deformation of the whole elastic domain.
Therefore, for all the test cases hereinafter, the $L^2$ {error norm} and $H^1$ {error semi-norm} of
the displacement obtained numerically are computed taking a more refined solution
as a reference. For this bidimensional test case, the mesh size
of the refined solution $h_{ref}$ is such that, for all the discretizations,
$4 h_{ref} \leq h$, where $h$ is the size of the mesh considered.
Additionally, as it is shown in Figure \ref{fig:mesh1_a}, the mesh is finer in the
vicinity of the potential contact zone. The knot vector values
are defined such that $80\%$ of the knot spans are located within $10\%$ of the total length of the knot vector.

In particular, the analysis of this example focuses on the effect of the
interpolation order on the quality of contact stress distribution.
Thus, in Figure \ref{fig:mesh1_b} we  compare the pressure reference solution.
{with the Lagrange multiplier values computed at the control points, \textit{i.e.} its constant values, and a post-processing which consists in a $P1$ re-interpolation.}
The dimensionless contact pressure $p/p_0$ is plotted respect to the normalized
coordinate $r/a$.
% with $r$ the distance respect to the centre of the contact area.
{The results are very good: the maximum pressure computed
and the pressure distribution, even across the boundary of the contact region (on the contact and non contact zones), are close to the analytical solution.}
% \begin{figure}[!ht]
%  \centering
%  \includegraphics[width= 8cm]{figure/2D_Q2_P0point003/ref_plot-eps-converted-to.pdf}
%  \caption{reference solution for $P = 0.003$.} \label{fig:ref sol 0point003}
%\end{figure} 

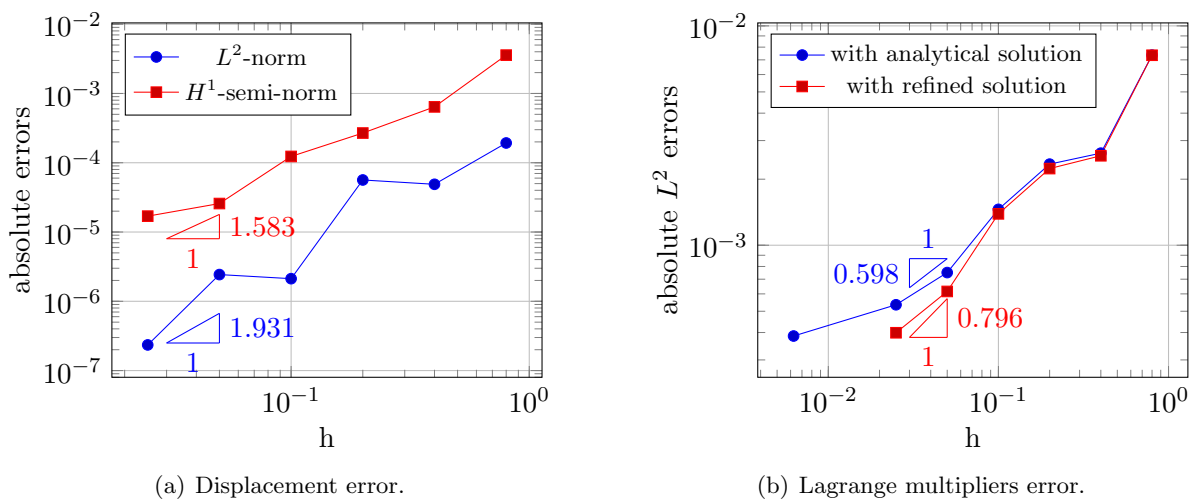
\begin{figure}[!ht]
 \centering
 \subfigure[Displacement error.]{\label{fig:disp 0point003_a}
\begin{tikzpicture}
\begin{loglogaxis}[width=7.3cm,xlabel={h},ymin=0.8e-07,ylabel={absolute errors},legend pos={north west},grid=major]

\addplot table[x=h,y=L2_abs] {figure_latex_figure_2d_hertz_p_0_003_N2_S0_data_disp_N2_S0_p_0point003.res};
\addplot table[x=h,y=H1_abs] {figure_latex_figure_2d_hertz_p_0_003_N2_S0_data_disp_N2_S0_p_0point003.res};

\slopeTriangleAbove{0.03}{0.05}{2.5e-07}{6.7039e-07}{1.931}{blue}; 
\slopeTriangleAbove{0.03}{0.05}{8e-06}{1.796e-05}{1.583}{red}; 
\legend{ \footnotesize $L^2$-norm, \footnotesize $H^1$-{semi-}norm}
\end{loglogaxis}
\end{tikzpicture}
} \hfill
 \subfigure[Lagrange multipliers error.]{\label{fig:disp 0point003_b}

\begin{tikzpicture}
\begin{loglogaxis}[width=7.3cm,xlabel={h},ymin=2.5e-04,ylabel={absolute $L^2$ errors},legend pos={north west},grid=major]

\addplot table[x=h_mult_ana,y=L2_mult_abs_ana] {figure_latex_figure_2d_hertz_p_0_003_N2_S0_data_mult_analytical_N2_S0_p_0point003.res};
\addplot table[x=h_mult_ref,y=L2_mult_abs_ref] {figure_latex_figure_2d_hertz_p_0_003_N2_S0_data_mult_refined_N2_S0_p_0point003.res};

\slopeTriangleBelow{0.03}{0.05}{6.4e-04}{8.6865e-04}{0.598}{blue}; 
\slopeTriangleAbove{0.03}{0.05}{3.8e-04}{5.70657e-04}{0.796}{red}; 
\legend{ \footnotesize with analytical solution,\footnotesize with refined solution}
\end{loglogaxis}
\end{tikzpicture}
}
 \caption{2D Hertz contact problem with $N_2 / S_0$ method for an applied pressure $P=0.003$.
 Absolute displacement errors in $L^2$-norm and $H^1$-{semi-}norm and Lagrange multipliers error in $L^2$-norm, respect to analytical and refined numerical solutions.} 
 \label{fig:disp 0point003}
\end{figure}

In Figure \ref{fig:disp 0point003_a}, absolute errors in $L^2${-norm} and $H^1${-semi-norm}
for the $N_2 / S_0$ choice are shown. As expected, optimal convergence is obtained
for the displacement error in the $H^1$-{semi-}norm: the convergence rate is close to the expected $3/2$ value.
Nevertheless, the $L^2$-norm of the displacement error presents suboptimal convergence (close to $2$), but
according to Aubin-Nitsche's lemma in the linear case, the expected convergence rate is {lower than} $5/2$.
On the other hand, in Figure \ref{fig:disp 0point003_b} the 
$L^2$-norm of the Lagrange multipliers error is presented, the expected convergence rate is $1$.
Whereas a convergence rate close to $0.6$ is achieved when {we compare the numerical solution and the Hertz's analytical solution}, and close to $0.8$ is achieved when {we compare the numerical solution and the refined numerical solution}.

As a second example, we present the same test case but with significantly higher pressure applied
$P = 0.01$. Under these load conditions, the contact area is wider ($a =  0.15223$) and the
contact pressure higher ($p_0 = 0.083641$).
It can be considered that the ratio $a/R$ no longer satisfies the hypotheses of Hertz's theory.

\begin{figure}[!ht]
  \centering
  \subfigure[Stress magnitude distribution for the undeformed mesh.]{\label{fig:mesh2_a} \includegraphics[width= 8.5cm]{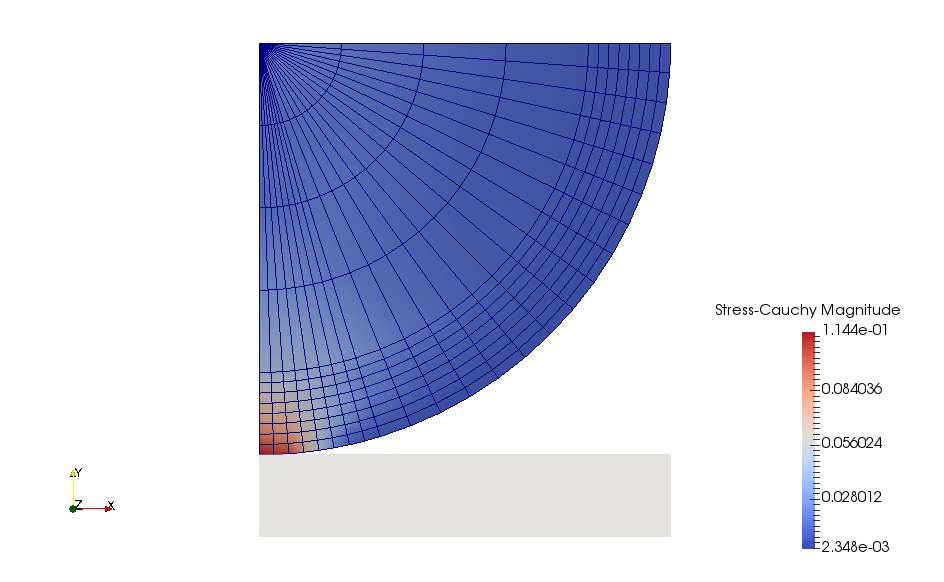}}\hfill
  \subfigure[Analytical and numerical contact pressure.]{\label{fig:mesh2_b} \includegraphics[width= 7cm]{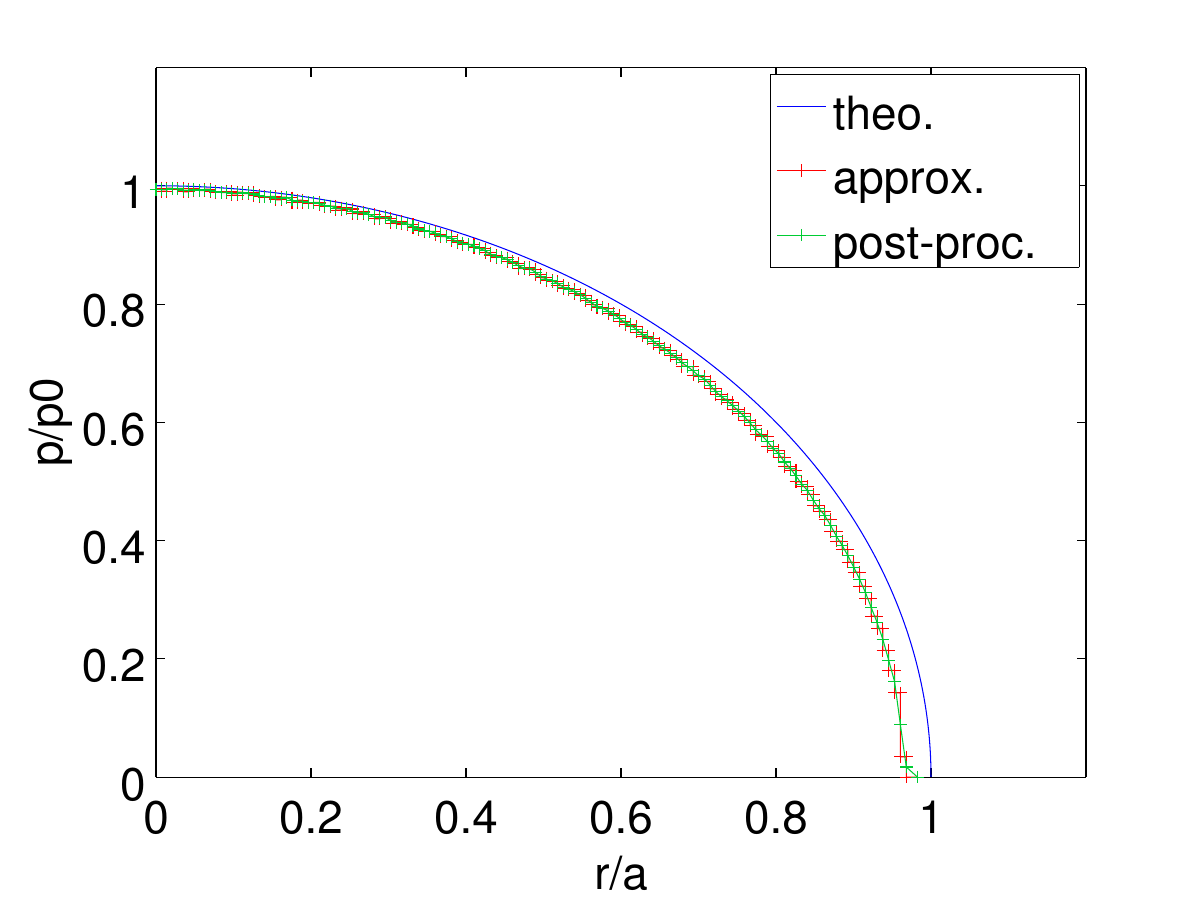}}
  \caption{2D Hertz contact problem with $N_2 / S_0$ method for a higher applied pressure ($P=0.01$).}
  \label{fig:mesh2}
\end{figure}

\begin{figure}[!ht]
 \centering
 \subfigure[Displacement error.]{\label{fig:disp 0point01_a} 
 
 \begin{tikzpicture}
\begin{loglogaxis}[width=7.3cm,ymin=1e-08,xmin=2e-02,xlabel={h},ylabel={absolute errors},legend pos={north west},grid=major]

\addplot table[x=h,y=L2_abs] {figure_latex_figure_2d_hertz_p_0_01_N2_S0_data_disp_N2_S0_p_0point01.res};
\addplot table[x=h,y=H1_abs] {figure_latex_figure_2d_hertz_p_0_01_N2_S0_data_disp_N2_S0_p_0point01.res};

\slopeTriangleAbove{0.03}{0.05}{7.5e-08}{2.52320e-07}{2.375}{blue}; 
\slopeTriangleAbove{0.03}{0.05}{1.7e-05}{3.6114e-05}{1.475}{red}; 
\legend{ \footnotesize $L^2$-norm, \footnotesize $H^1$-{semi-}norm}

\end{loglogaxis}
\end{tikzpicture}
 
 }\hfill
 \subfigure[Lagrange multipliers error.]{\label{fig:disp 0point01_b}
 
\begin{tikzpicture}
\begin{loglogaxis}[every axis y label/.style={at={(ticklabel cs:0.5)},rotate=90,anchor=center},width=7.3cm,xlabel={h},ymin=4e-04,ylabel={absolute $L^2$ errors},legend pos={north west},xmin=1.5e-3,grid=major]

\addplot table[x=h_mult_ana,y=L2_mult_abs_ana] {figure_latex_figure_2d_hertz_p_0_01_N2_S0_data_mult_analytical_N2_S0_p_0point01.res};
\addplot table[x=h_mult_ref,y=L2_mult_abs_ref] {figure_latex_figure_2d_hertz_p_0_01_N2_S0_data_mult_refined_N2_S0_p_0point01.res};

\slopeTriangleBelow{0.0062500}{0.012}{0.0011}{0.00113498767}{0.048}{blue}; 
\slopeTriangleAbove{0.035}{0.05}{5.6e-04}{7.8083e-04}{0.932}{red}; 
\legend{ \footnotesize with analytical solution,\footnotesize with refined solution}

\end{loglogaxis}
\end{tikzpicture}  
  
  }
\caption{2D Hertz contact problem with $N_2 / S_0$ method for an applied pressure $P=0.01$.
 Absolute displacement errors in $L^2$-norm and $H^1$-{semi-}norm and Lagrange multipliers error in $L^2$-norm, respect to analytical and refined numerical solutions.} 
 \label{fig:disp 0point01}
\end{figure}
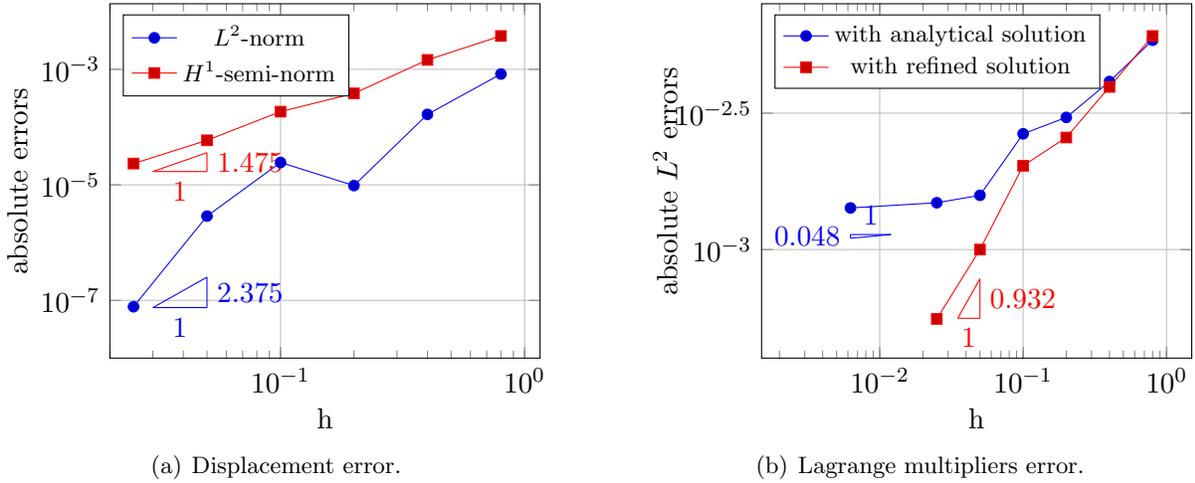

In the same way as before, Figure \ref{fig:mesh2} shows the stress tensor magnitude
and computed contact pressure.
Figure \ref{fig:disp 0point01_a} shows the displacement absolute error in $L^2$-norm and $H^1$-{semi-}norm
for $N_2 / S_0$ method. As expected, optimal convergence is obtained in the $H^1$-{semi-}norm,
(the convergence rate is close to $1.5$) and, {while, for the $L^2$-norm we obtain a better rate (as expected by the Aubin-Nische's lemma) which can hardly be estimated precisely from the graph.}
On the other hand, in Figure \ref{fig:disp 0point01_b} it can be seen that the $L^2$-norm of the error of the Lagrange multipliers
evidences a suboptimal behaviour: the error, that initially decreases, remains constant for smaller values of $h$.
It may due to the choice of an excessively large normal pressure:
the approximated solution converges, but not to the analytical solution, that is no longer valid.
Indeed, when compared to a refined numerical solution (Figure \ref{fig:disp 0point01_b}),
the computed Lagrange multipliers solution converges optimally. As it was pointed out above,
for these examples the displacement solution error is computed respect to a more refined numerical solution, therefore, this effect
is not present in displacement results.

\subsection{Three-dimensional Hertz problem}
\label{subsec:Hertz problems 3d}
\begin{figure}[!ht]
  \centering
  \subfigure[magnitude for the undeformed mesh.]{\label{fig:3d p=0.0001_a} \includegraphics[width= 8.5cm]{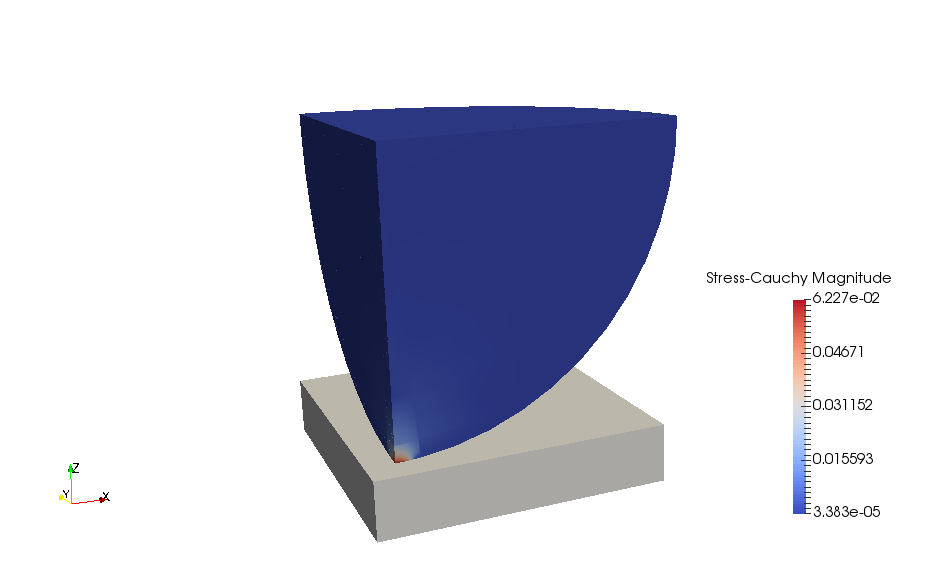}}\hfill
  \subfigure[Analytical and numerical contact pressure for $h=0.1$.]{\label{fig:3d p=0.0001_b} \includegraphics[width= 7cm]{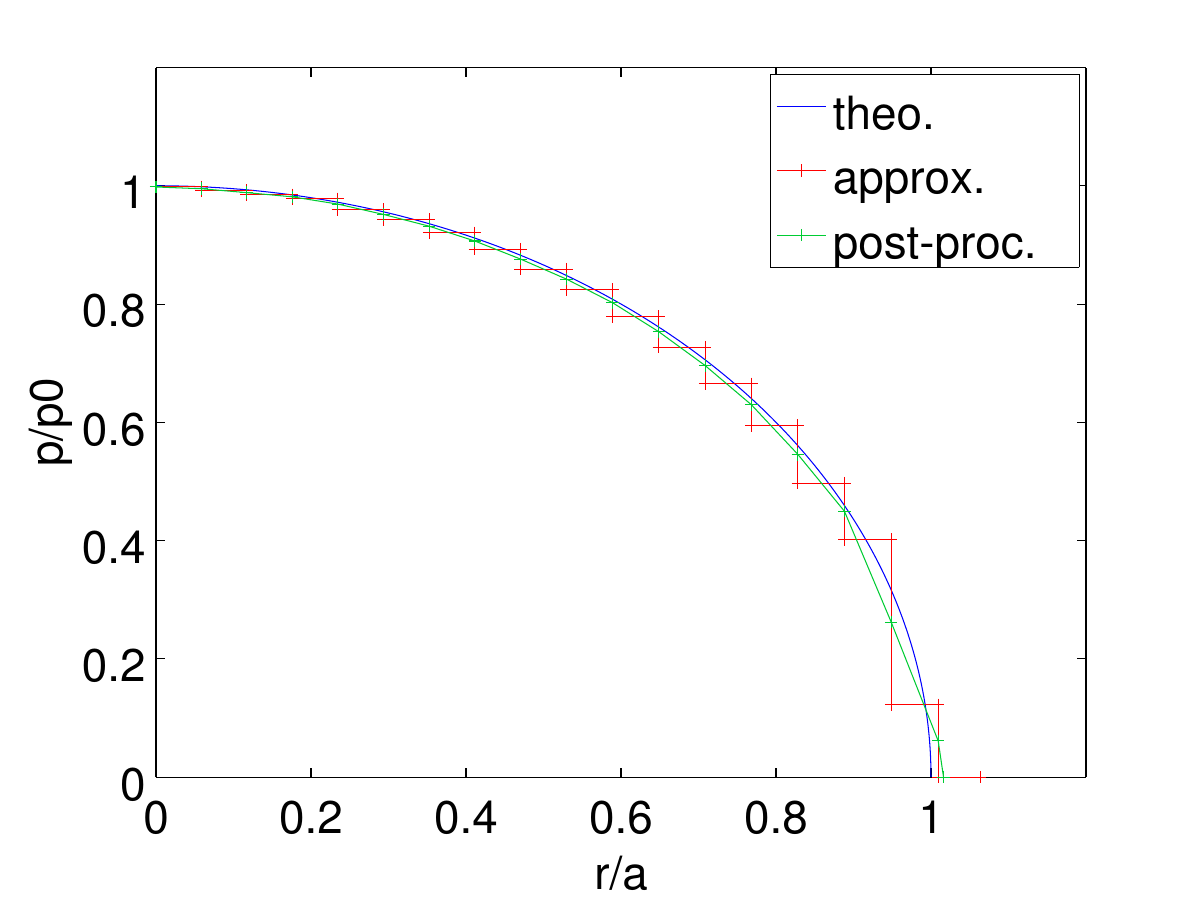}}
  \caption{3D Hertz contact problem with $N_2 / S_0$ method for an applied pressure $P=10^{-4}$.}
  \label{fig:3d p=0.0001}
\end{figure}

In this section, the three-dimensional frictionless Hertz problem is studied.
It consists in a hemispherical elastic body with radius $R$ that contacts against a horizontal rigid plane as a consequence of
an uniform pressure $P$ applied on the top face (see Figure \ref{fig:3d p=0.0001_a}).
Hertz's theory predicts that the contact region is a circle of radius $a = (3R^3P(1- \nu^2)/4 E)^{1/3}$
and the contact pressure follows a hemispherical distribution
$p(r)=p_0 \sqrt{1-r^2/a^2}$, with $p_0 = 3R^2P/ 2 a^2$, being $r$ the distance to the centre of the circle
(see\cite{Hertz1882,johnsonbook}).
In this case, for the chosen values $R=1$, $E=1$, $\nu=0.3$ and $P=10^{-4}$, the contact radius is $a=0.059853$ and the maximum pressure $p_0 = 0.041872$.
As in the two-dimensional case, Hertz's theory relies on the hypothesis that $a$ is small compared to $R$ and the deformations are small.

Considering the problem axial symmetry, the test is reproduced using an octant of sphere with proper boundary conditions.
Figure \ref{fig:3d p=0.0001_a} shows the problem setup and the magnitude of the computed stresses.
As in the 2D case, in order to achieve more accurate results in the contact region, the mesh is refined in the vicinity of the potential contact zone.
The knot vectors are defined such as $75\%$ of the elements are located within $10\%$ of the total length of the knot vector.

\begin{figure}[!ht]
  \centering
  \subfigure[$h=0.4$.]{\label{fig:3d cp 0.0001_a} \includegraphics[width= 7cm]{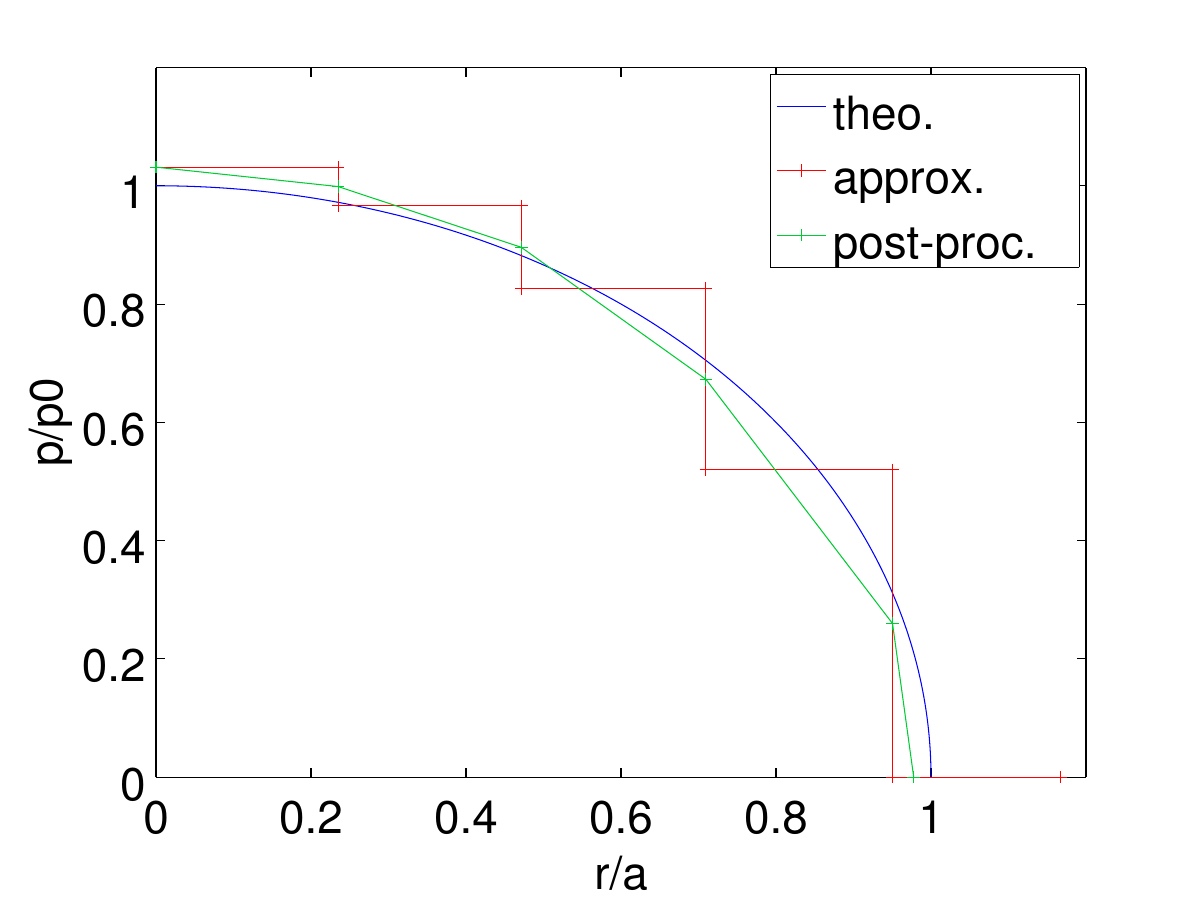}}\hfill
  \subfigure[$h=0.2$.]{\label{fig:3d cp 0.0001_b} \includegraphics[width= 7cm]{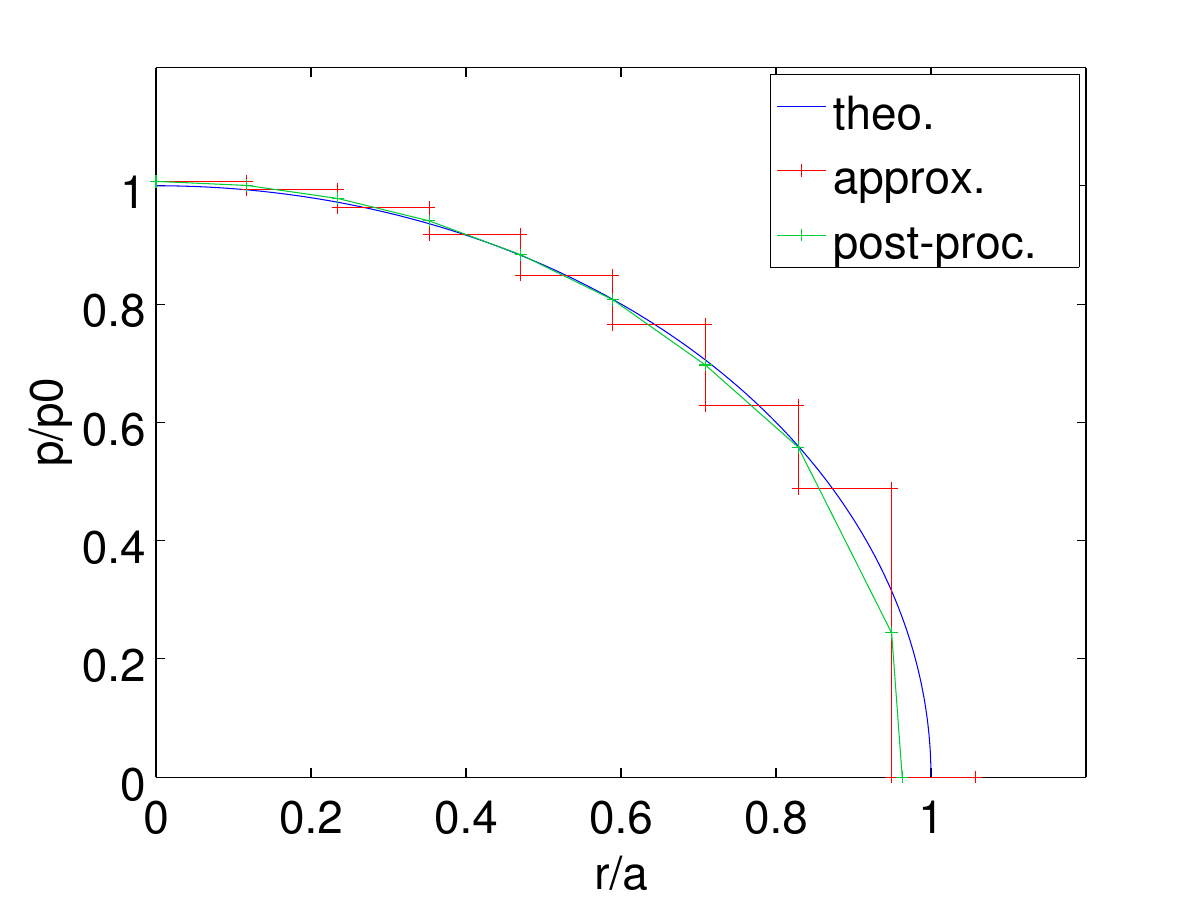}}
  \caption{3D Hertz contact problem with $N_2 / S_0$ method for an applied pressure $P=10^{-4}$. Contact pressure solution at control points.}
  \label{fig:3d cp 0.0001}
\end{figure}

{In Figure \ref{fig:3d p=0.0001_b}, we compare the Hertz's solution with the computed contact pressure at control points
and a $P1$ re-interpolation of those values, for a mesh with size $h=0.1$.}
On the other hand, in Figure \ref{fig:3d cp 0.0001} the contact pressure is shown at control points for mesh sizes $h=0.4$ and $h=0.2$. As it can be
appreciated, good agreement between the analytical and computed pressure is obtained in all cases.

As in the previous test, the displacement of the deformed elastic body is not fully described by the Hertz's theory.
Therefore, the $L^2$ error norm and $H^1$ error {semi}-norm of the displacement are evaluated
by comparing the obtained solution with a finer refined case. Nonetheless,
Lagrange multipliers computed solutions are compared with the analytical contact pressure.
In this test case, the size of the refined mesh is $h_{ref}=0.1175$ ($0.0025$ in the contact region), and
it is such as $2 h_{ref} \leq h$.

\begin{figure}[!ht]
  \centering
  \subfigure[Displacement error.]{\label{fig:3d  error disp 1_a}
  
  \begin{tikzpicture}
\begin{semilogyaxis}[width=7.3cm,ymin=0.9e-04,xmin=2e-02,xlabel={h},ylabel={absolute errors},legend pos={north west},grid=major]

\addplot table[x=h,y=L2_abs] {figure_latex_figure_3d_hertz_p_0_0001_N2_S0_data_disp_3D_N2_S0_p_0point0001.res};
\addplot table[x=h,y=H1_abs] {figure_latex_figure_3d_hertz_p_0_0001_N2_S0_data_disp_3D_N2_S0_p_0point0001.res};

\slopeTriangleAbove{0.3}{0.4}{1.14e-04}{1.5806e-04}{1.26}{blue}; 
\slopeTriangleBelow{0.3}{0.4}{0.0003}{0.000347}{0.5}{red}; 
\legend{ \footnotesize $L^2$-norm, \footnotesize $H^1$-{semi-}norm %\pgfmathprintnumber{\pgfplotstableregressiona},
}

\end{semilogyaxis}
\end{tikzpicture}
  
  }\hfill
  \subfigure[Lagrange multipliers error.]{\label{fig:3d  error disp 1_b}
  
\begin{tikzpicture}
\begin{semilogyaxis}[every axis y label/.style={at={(ticklabel cs:0.5)},rotate=90,anchor=center},ymax=2e-03,width=7.3cm,xlabel={h},ylabel={absolute $L^2$ errors},legend pos={north west},xmin=1.5e-3,grid=major]

\addplot table[x=h_mult_ana,y=L2_mult_abs_ana] {figure_latex_figure_3d_hertz_p_0_0001_N2_S0_data_mult_analytical_3D_N2_S0_p_0point0001.res};
\addplot table[x=h_mult_ref,y=L2_mult_abs_ref] {figure_latex_figure_3d_hertz_p_0_0001_N2_S0_data_mult_refined_3D_N2_S0_p_0point0001.res};

\slopeTriangleBelow{0.2500}{0.3500}{6.7e-04}{8.34e-04}{0.649}{blue}; 
\slopeTriangleAbove{0.2500}{0.3500}{3.9e-04}{5.11e-04}{0.806}{red}; 
\legend{ \footnotesize with analytical solution,\footnotesize with refined solution}

\end{semilogyaxis}
\end{tikzpicture}  
  
}
  \caption{3D Hertz contact problem with $N_2 / S_0$ method for an applied pressure $P=10^{-4}$.
   Absolute displacement errors in $L^2$-norm and $H^1$-{semi-}norm and Lagrange multipliers error in $L^2$-norm, respect to analytical and refined numerical solutions.} 
  \label{fig:3d  error disp 1}
\end{figure}
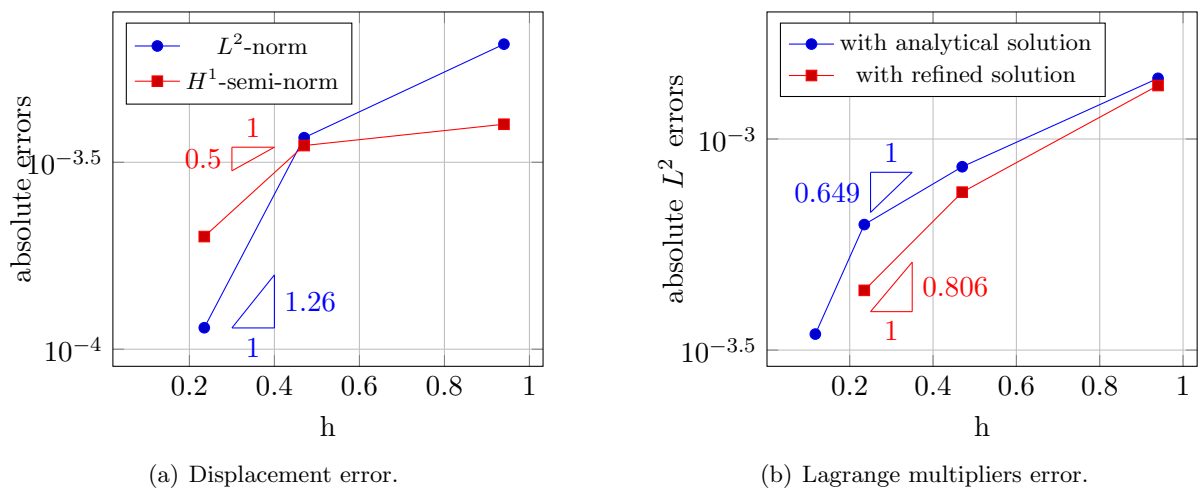

In Figure \ref{fig:3d  error disp 1_a} the displacement error norms are reported. As it can {be} seen, they present suboptimal convergence rates
both in the $L^2$-norm and $H^1$-{semi-}norm. Convergence rates are close to $1.26$ and $0.5$, respectively.
The large mesh size of the numerical reference solution $h_{ref}$, limited by our computational resources, seems
to be the cause of these suboptimal results. {Due to the coarse reference mesh, the presented rates are only pre-asymptotic.}
Better behaviour is observed for the Lagrange multipliers error (Figure \ref{fig:3d  error disp 1_b}).

\begin{figure}[!ht]
  \centering
  \subfigure[Stress magnitude distribution for the undeformed mesh.]{\label{fig:3d p = 0.0005_a} \includegraphics[width= 8.5cm]{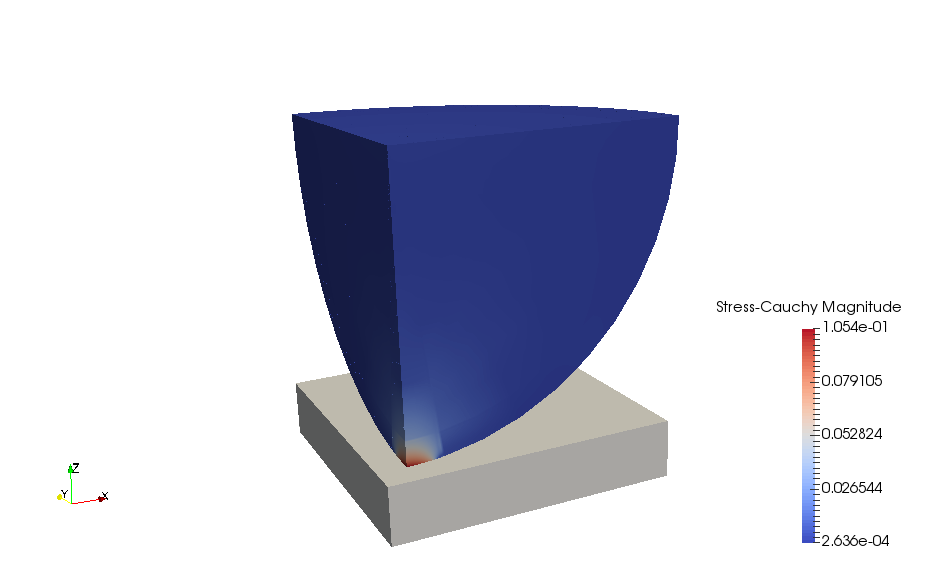}}\hfill
  \subfigure[Analytical and numerical contact pressure for $h = 0.1$.]{\label{fig:3d p = 0.0005_b} \includegraphics[width= 7cm]{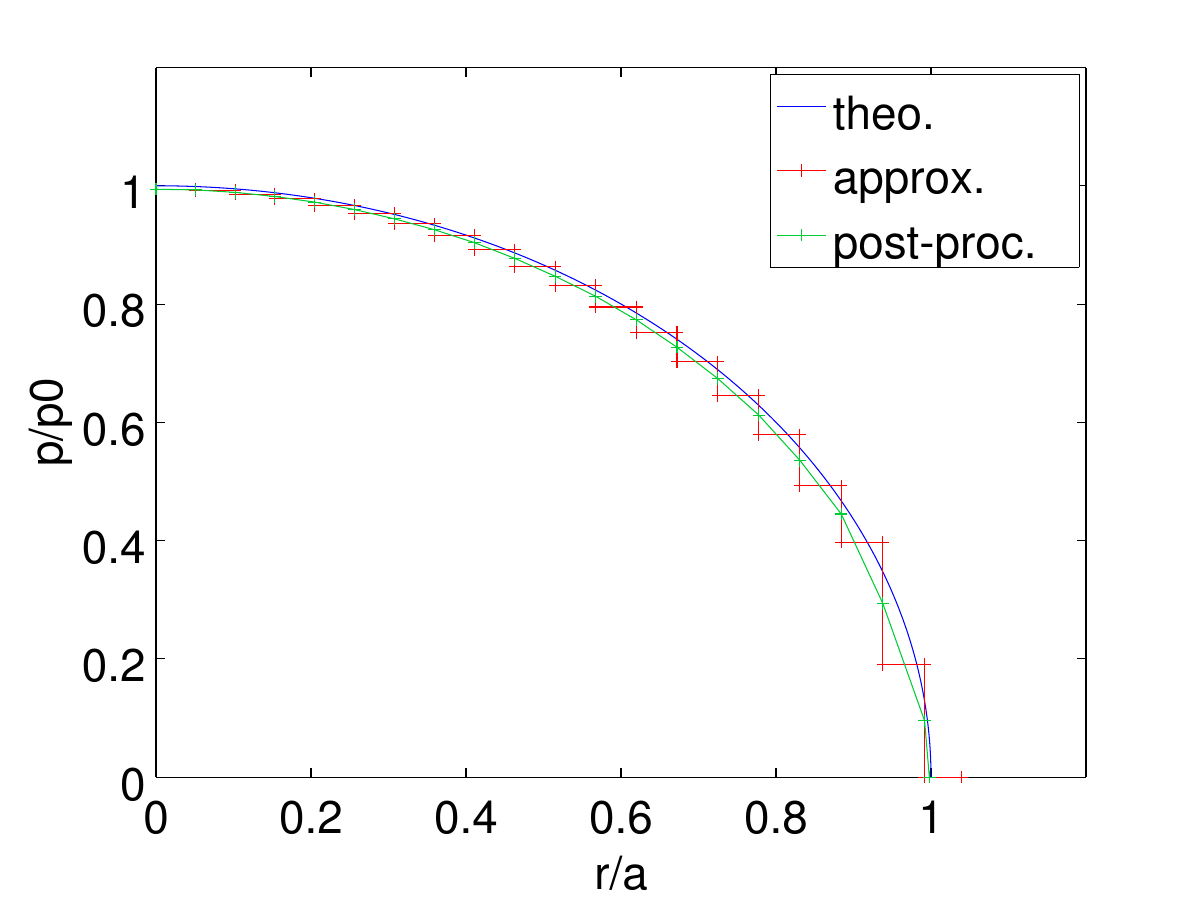}}
  \caption{3D Hertz contact problem with $N_2 / S_0$ method for a higher pressure ($P=5\cdot10^{-4}$).}
  \label{fig:3d p = 0.0005}
\end{figure}

\begin{figure}[!ht]
  \centering
  \subfigure[$h=0.4$.]{\label{fig:3d cp 0.0005_a} \includegraphics[width= 7cm]{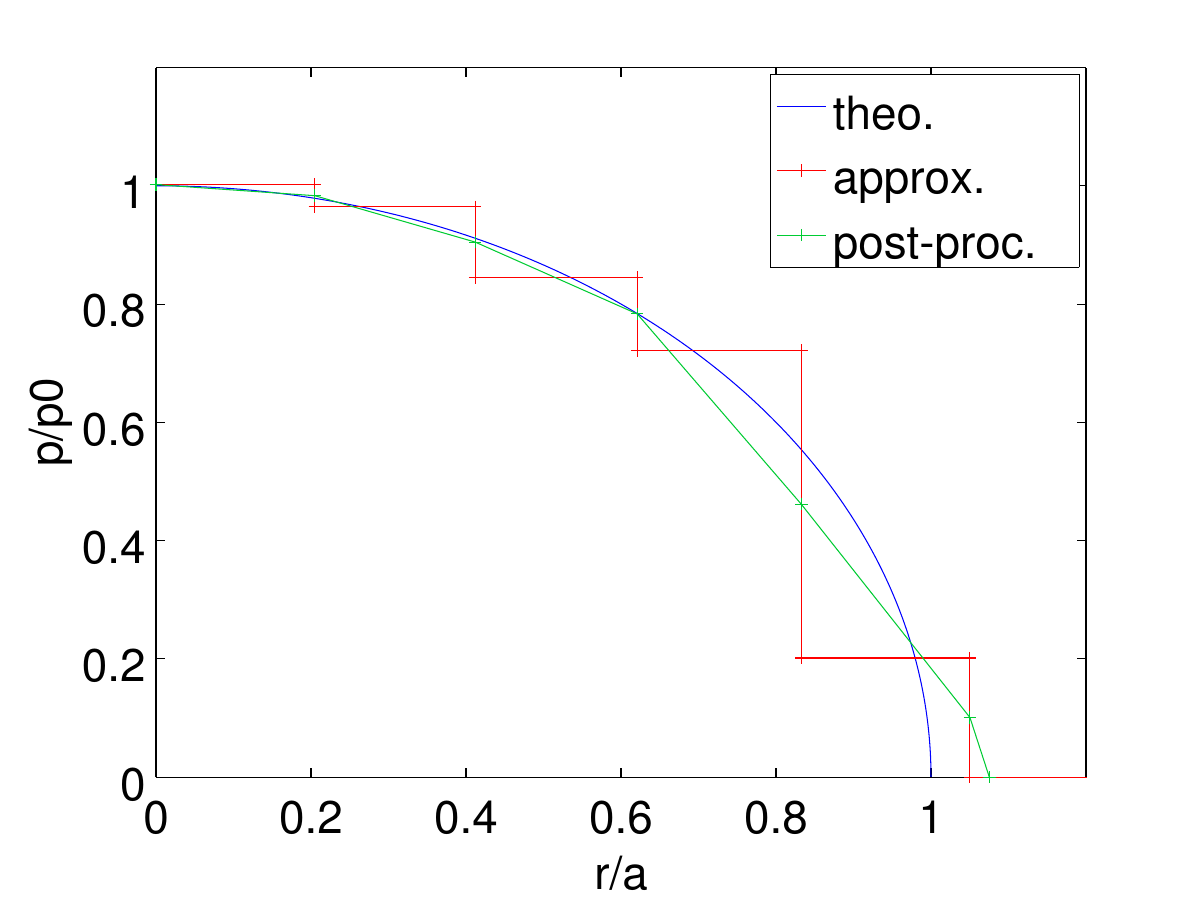}}\hfill
  \subfigure[$h=0.2$.]{\label{fig:3d cp 0.0005_b} \includegraphics[width= 7cm]{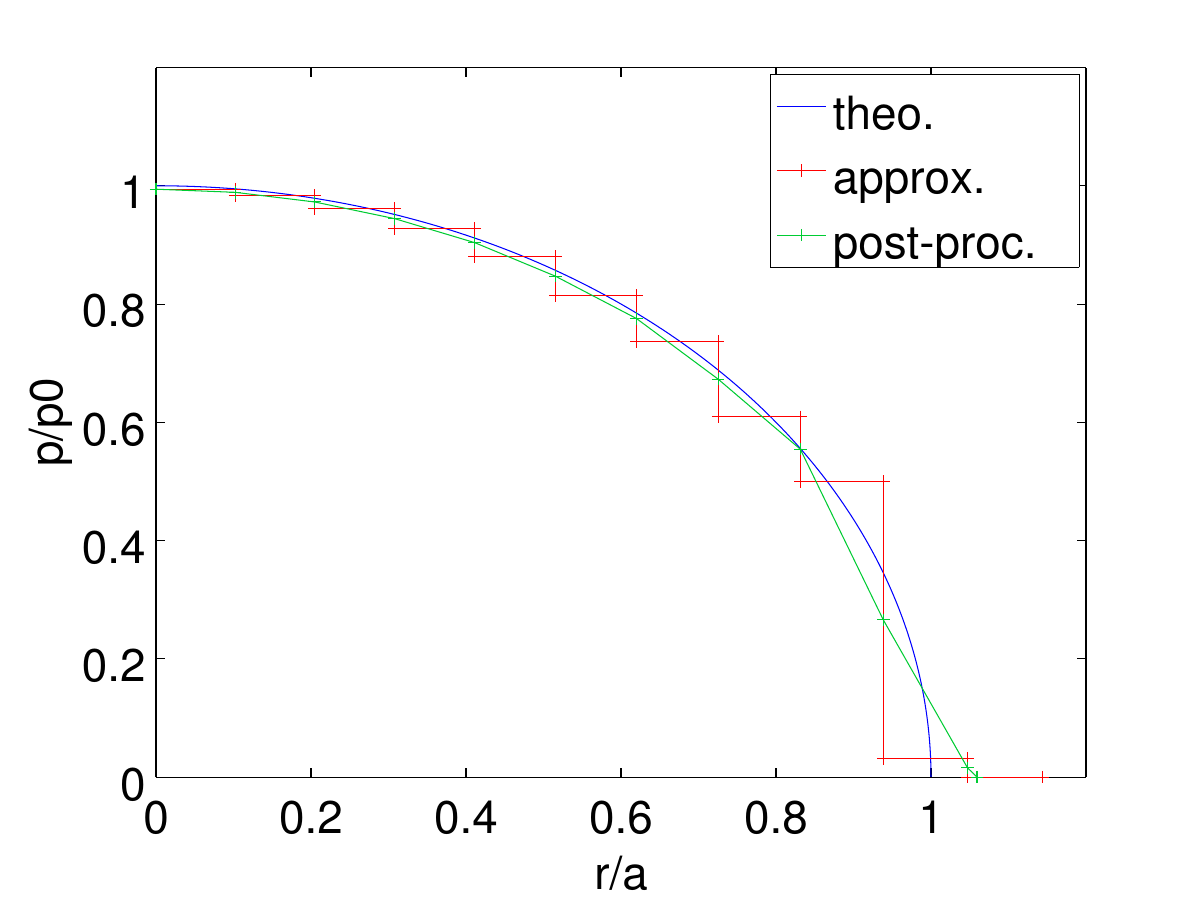}}
  \caption{3D Hertz contact problem with $N_2 / S_0$ method for an applied pressure $P=5\cdot10^{-4}$. Contact pressure solution at control points.}
  \label{fig:3d cp 0.0005}
\end{figure}

By considering a higher pressure ($P=5\cdot10^{-4}$), the radius of the contact
zone becomes larger ($a = 0.10235$), and the ratio $a/R$ does not satisfy the theory hypotheses.
Figure \ref{fig:3d p = 0.0005} shows the stress magnitude and contact pressure at the control points
for a given mesh.
Similarly, in Figure \ref{fig:3d cp 0.0005} the analytical contact pressure is compared
with the computed Lagrange multipliers values associated to the control points for different meshes. Satisfactory results
are observed in all cases.

As in the previous test, the coarse value of the reference mesh size $h_{ref}$ seems
to be the cause of the suboptimal convergence of the displacement shown in Figure \ref{fig:3d  error disp 2_a}.
An optimal convergence is observed for the Lagrange multipliers error in the $L^2$-norm (see Figure \ref{fig:3d  error disp 2_b}). However, due to the coarse value of the mesh size, we do not observe the expected threshold of the $L^2$ error for the Lagrange multipliers between the analytical and approximate solutions.
 \begin{figure}[!ht]
  \centering
  \subfigure[Displacement error.]{\label{fig:3d  error disp 2_a} 
  
\begin{tikzpicture}
\begin{semilogyaxis}[width=7.3cm,ymin=0.9e-04,ymax=5e-03,xmin=2e-02,xlabel={h},ylabel={absolute errors},legend pos={north west},grid=major]

\addplot table[x=h,y=L2_abs] {figure_latex_figure_3d_hertz_p_0_0005_N2_S0_data_disp_3D_N2_S0_p_0point0005.res};
\addplot table[x=h,y=H1_abs] {figure_latex_figure_3d_hertz_p_0_0005_N2_S0_data_disp_3D_N2_S0_p_0point0005.res};

\slopeTriangleAbove{0.3}{0.4}{0.000188257}{0.000285}{1.441}{blue}; 
\slopeTriangleBelow{0.3}{0.4}{0.00095}{0.00115}{0.659}{red}; 
\legend{ \footnotesize $L^2$-norm, \footnotesize $H^1$-{semi-}norm %\pgfmathprintnumber{\pgfplotstableregressiona},
}

\end{semilogyaxis}
\end{tikzpicture}  
  }\hfill
  \subfigure[Lagrange multipliers error.]{\label{fig:3d  error disp 2_b}

\begin{tikzpicture}
\begin{semilogyaxis}[every axis y label/.style={at={(ticklabel cs:0.5)},rotate=90,anchor=center},width=7.3cm,xlabel={h},ymax=0.8e-02,ylabel={absolute $L^2$ errors},legend pos={north west},xmin=1.5e-3,grid=major]

\addplot table[x=h_mult_ana,y=L2_mult_abs_ana] {figure_latex_figure_3d_hertz_p_0_0005_N2_S0_data_mult_analytical_3D_N2_S0_p_0point0005.res};
\addplot table[x=h_mult_ref,y=L2_mult_abs_ref] {figure_latex_figure_3d_hertz_p_0_0005_N2_S0_data_mult_refined_3D_N2_S0_p_0point0005.res};

\slopeTriangleBelow{0.2500}{0.3500}{1.7659e-03}{2.4220e-03}{0.939}{blue}; 
\slopeTriangleAbove{0.2500}{0.3500}{0.0011023}{0.001558}{1.029}{red}; 
\legend{ \footnotesize with analytical solution,\footnotesize with refined solution}

\end{semilogyaxis}
\end{tikzpicture}

   }
  \caption{3D Hertz contact problem with $N_2 / S_0$ method for an applied pressure $P=5\cdot10^{-4}$.
   Absolute displacement errors in $L^2$-norm and $H^1$-{semi-}norm and Lagrange multipliers error in $L^2$-norm, respect to analytical and refined numerical solutions.} 
 \label{fig:3d  error disp 2}
\end{figure}

\subsection{Two-dimensional Hertz problem with large deformations} \label{subsec:large deformation}

Finally, in this section the two-dimensional frictionless Hertz problem is studied considering large deformations and strains.
For that purpose, a Neo-Hookean material constitutive law {(an hyper-elastic law that considers finite strains)} with Young's modulus $E = 1$ and Poisson's ratio $\nu = 0.3$,
has been used for the deformable body.

%\emph{\pablocor{We should address here the comment of Reviewer 1 about linear/non-linear contact conditions. I don't see his/her point referring to the Equation (2) for the linear contact conditions. In my view, both the small and large strain problems use non-linear contact conditions in the sense that that the contact state is computed referred to the deformed geometry of the body, isn't it? Anyway, I think he/she is right and is worth mention it explicitly here. What do you think?}}

\begin{figure}[!ht]
  \centering
  \subfigure[Stress magnitude distribution for the deformed mesh.]{\label{fig:mesh3_a} \includegraphics[width= 8.5cm]{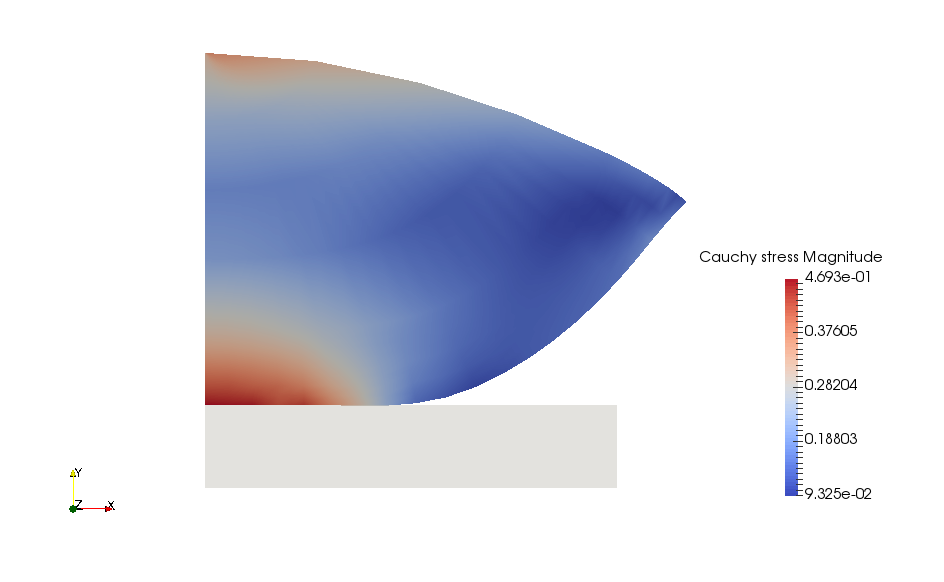}}\hfill
  \subfigure[Analytical and numerical contact pressure for $h = 0.1$.]{\label{fig:mesh3_b} \includegraphics[width= 7cm]{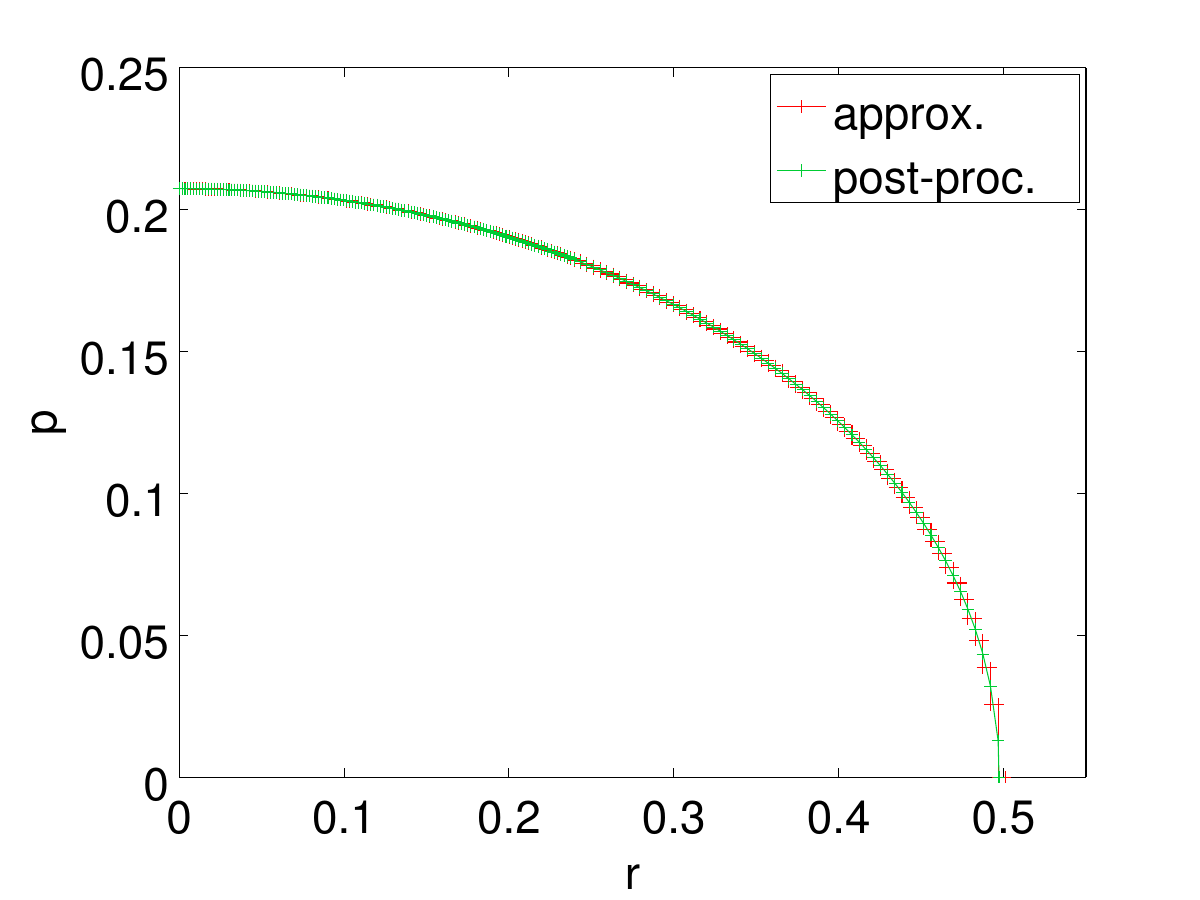}}
  \caption{2D large deformation Hertz contact problem with $N_2 / S_0$ method for an applied pressure $P = 0.1$.}
  \label{fig:mesh3}
\end{figure}

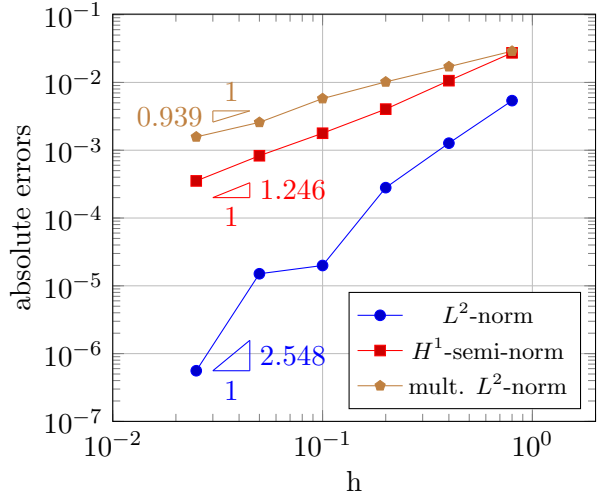
\begin{figure}[!ht]
  \centering

\begin{tikzpicture}
\begin{loglogaxis}[width=8cm, height={7cm} ,ymin=1e-07,xmin=1e-02,xmax=2e0,
xlabel={h},ylabel={absolute errors},legend pos={south east},grid=major]

\addplot table[x=h,y=L2_abs] {figure_latex_figure_data_large_neumann_data_large_neumann.res};
\addplot table[x=h,y=H1_abs] {figure_latex_figure_data_large_neumann_data_large_neumann.res};
\addplot[brown,mark=pentagon*] table[x=h,y=L2_mult_abs] {figure_latex_figure_data_large_neumann_data_large_neumann.res};

\slopeTriangleAbove{0.03}{0.045}{5.58631e-07}{1.5697e-06}{2.548}{blue}; 
\slopeTriangleAbove{0.03}{0.045}{0.0002}{0.0003315}{1.246}{red}; 
\slopeTriangleBelow{0.03}{0.045}{0.0026}{0.003805}{0.939}{brown}; 

\legend{ \footnotesize $L^2$-norm, \footnotesize $H^1$-{semi-}norm,\footnotesize mult. $L^2$-norm}
\end{loglogaxis}
\end{tikzpicture}

  \caption{2D large deformation Hertz contact problem with $N_2 / S_0$ method for an applied pressure $P=0.1$.
   Absolute displacement errors in $L^2$-norm and $H^1$-{semi-}norm and Lagrange multipliers error in $L^2$-norm.} 
  \label{fig:disp 0point1}
\end{figure}

As in Section \ref{subsec:Hertz problems}, the performance
of the $N_2 / S_0$ method is analysed and the problem is modelled as an elastic quarter of disc with proper boundary conditions.
The considerations made about the mesh size in Section \ref{subsec:Hertz problems} are also valid for the present case.
The radius of the cylinder is $R=1$ and the applied pressure $P=0.1$ (ten times higher than the one considered in Section \ref{subsec:Hertz problems}).
In this large deformation framework the exact solution is unknown: the error of the computed displacement and Lagrange multipliers are studied
taking a refined numerical solution as reference.
Figure \ref{fig:mesh3} shows the final deformation of the elastic body and the computed contact pressure.
In Figure \ref{fig:disp 0point1}, the displacement and multiplier errors are reported. It can be seen that the obtained displacement presents
optimal convergence both in $L^2$-norm and $H^1$-{semi-}norm; analogously, optimal convergence is also achieved
for the computed Lagrange multipliers.

\begin{figure}[!ht]
  \centering
  \subfigure[Stress magnitude distribution for the deformed mesh.]{\label{fig:mesh4_b} \includegraphics[width= 8.5cm]{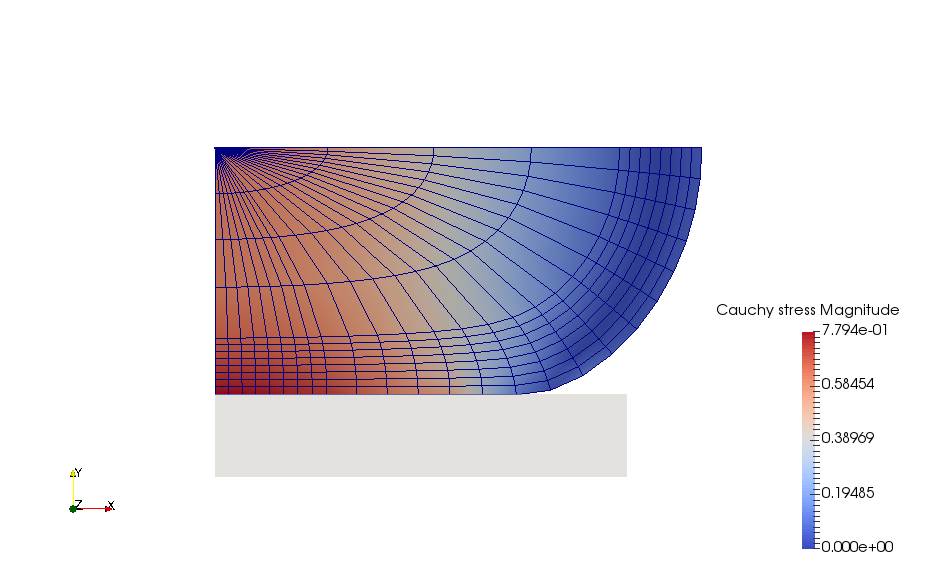}}\hfill
  \subfigure[Analytical and numerical contact pressure for $h = 0.1$]{\label{fig:mesh4_b} \includegraphics[width= 7cm]{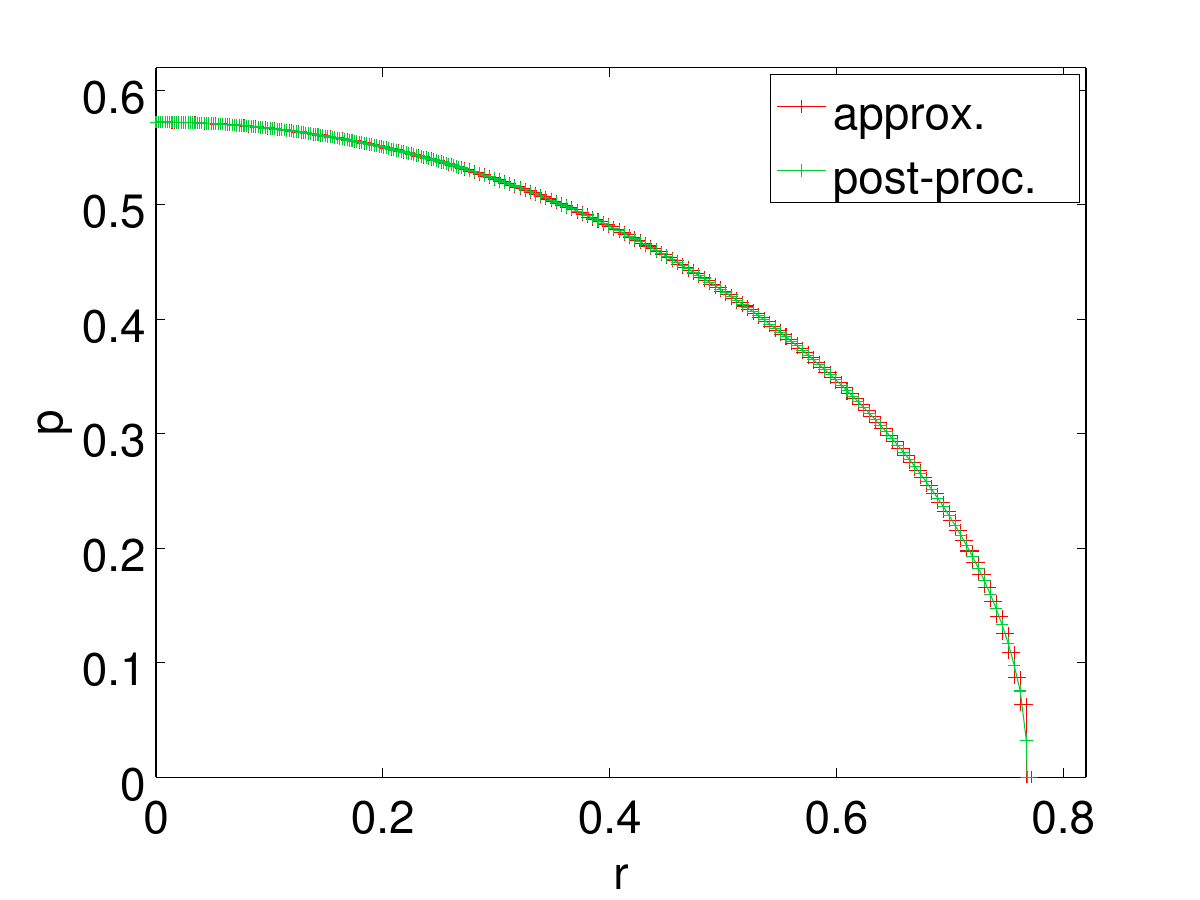}}
  \caption{2D large deformation Hertz contact problem with $N_2 / S_0$ method with a uniform downward displacement $u_y = -0.4$.}
  \label{fig:mesh4}
\end{figure}

As a last example, the same large deformation Hertz problem is considered, but modifying its boundary conditions:
instead of pressure, a uniform downward displacement $u_y=-0.4$ is applied on the top surface of the cylinder.
The large deformation of the body and computed contact pressure are presented in Figure \ref{fig:mesh4}.
As in the previous case (large deformation with applied pressure) optimal results are obtained for
the computed displacement and Lagrange multipliers.
\begin{figure}[!ht]
  \centering

\begin{tikzpicture}
\begin{loglogaxis}[width=8cm, height={7cm} ,ymin=4e-07,xmin=1.2e-02,xmax=2e0,
xlabel={h},ylabel={absolute errors},legend pos={south east},grid=major]

\addplot table[x=h,y=L2_abs] {figure_latex_figure_data_large_dirichlet_data_large_dirichlet.res};
\addplot table[x=h,y=H1_abs] {figure_latex_figure_data_large_dirichlet_data_large_dirichlet.res};
\addplot[brown,mark=pentagon*] table[x=h,y=L2_mult_abs] {figure_latex_figure_data_large_dirichlet_data_large_dirichlet.res};

\slopeTriangleAbove{0.03}{0.045}{1.4e-06}{3.711e-06}{2.404}{blue}; 
\slopeTriangleAbove{0.03}{0.045}{0.0002}{0.000364}{1.476}{red}; 
\slopeTriangleBelow{0.03}{0.045}{0.007}{0.0101}{0.908}{brown}; 

\legend{ \footnotesize $L^2$-norm, \footnotesize $H^1$-{semi-}norm,\footnotesize mult. $L^2$-norm}
\end{loglogaxis}
\end{tikzpicture}

  \caption{2D large deformation Hertz contact problem with $N_2 / S_0$ method with a uniform downward displacement $u_y = -0.4$.
   Absolute displacement errors in $L^2$-norm and $H^1$-{semi-}norm and Lagrange multipliers error in $L^2$-norm.} 
  \label{fig:disp 0point4}
\end{figure}
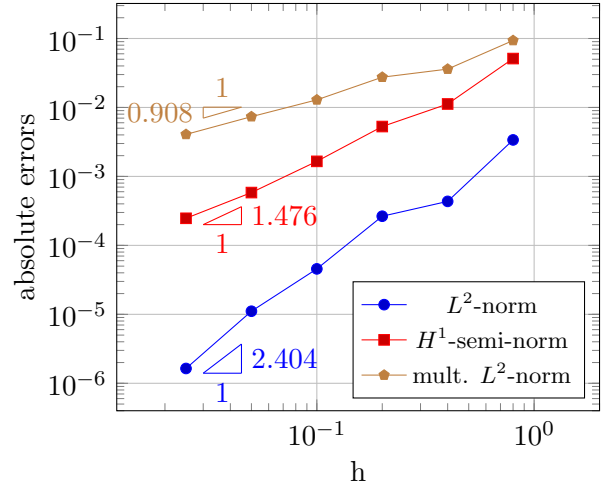\\

%%%%%%%%%%%%%%%%%%%%%%%%%%%%%%%%%%%%%%%%%%%%%%%%

\section*{Conclusions}
\label{sec: conclusion}
\addcontentsline{toc}{section}{Conclusions}

In this work, we present an {optimal} \textit{a priori} error estimate of unilateral contact problem frictionless between deformable body and rigid one. %we develop a Lagrange formulation for frictionless contact problem between rigid and deformable bodies for small and large deformation. This method combine IGA framework and active-set strategy between Lagrange multipliers and "local" gap. The contact status, for each parts of the workpiece, is updated at each iterations in order to expect converge. This method proposed by Alart and Curnier, which is combined with a Newton solution scheme, is applied to solve saddle point problem.

For the numerical point of view, we observe an {optimality} of this method for both variables, the displacement and the Lagrange multipliers. In our experiments, we use a NURBS of degree $2$ for the primal space and B-Spline of degree $0$ for the dual space. {Thanks to this choice of approximation spaces, we observe a stability of the Lagrange multipliers and a well approximation of the pressure in two-dimensional case and we observe a {sub-optimality} in three-dimensional case}. The {sub-optimality} observed in three-dimensional case may be due to the coarse mesh used. This NURBS based contact formulation seems to provide a robust description of large deformation.   

%Moreover, the use of a NURBS of degree $3$ for the primal space and B-Spline of degree $1$ for the dual space with multiplier method, active-set strategy and Newton's method does not converge in small and large deformations for the problems studied. %So we experiment augmented Lagrange method for this case. It seems to give good converge rates and stability.

%However, for optimization point of view between generalized objective functions under inequalities constrains, as large deformation contact problem, it seems to be interesting the choice of augmented Lagrange multipliers methods. Indeed, in some case the contact status can oscillate, which involves the non-convergence of the Newton scheme, or the Lagrangian method can converge to a "bad" solution. 

\section*{Acknowledgements}
\label{sec: acknowledgements}
\addcontentsline{toc}{section}{acknowledgements}
This work has been partially supported by Michelin under the contract A10-4087. A. Buffa acknowledges the support of the ERC Advanced grant no. 694515 and of the PRIN-MIUR project "Metodologie innovative nella modellistica differenziale numerica".

%\newpage
%\strut
%\newpage
%
%\newpage
%\strut
\newpage

 %Add in table of contents
\addcontentsline{toc}{section}{Appendices}
\section*{{Appendix 1.}}
\label{sec:appendix1}

In this appendix, we provide the ingredients needed to fully discretise the problem \eqref{eq:discrete_mixed_form} as well as its large deformation version that we have used in Section \ref{sec:sec3}. First we introduce the contact status, an active-set strategy for the discrete problem, and then the fully discrete problem. For the purpose of this appendix, we take notations suitable to large deformation and denote by $g_n$ the distance between the rigid and the deformable body. In small deformation, it holds $g_n(u)=u \cdot n$.
\subsection*{Contact status}
\label{subsec:contact stat}
Let us first deal with the contact status. The active-set strategy is defined in \cite{hueber-Wohl-05,hueber-Stadler-Wohlmuth-08} and is updated at each iteration of Newton. Due to the deformation, parts of the workpiece may come into contact or conversely may loose contact. This change of contact status changes the loading that is applied on the boundary of the mesh. This method is used to track the location of contact during the change in boundary conditions.\\
Let $K$ be a control point of the B-Spline space \eqref{def:space_mult}, let $(\Pi^h_\lambda \cdot )_K$ be the local projection defined in \eqref{def:l2proj_K} and let $P{\{}\lambda_K, (\Pi_\lambda^h g_n)_K{\}}$ be he operator defined component wise by:
\begin{itemize}
\item $\lambda_K = 0$, \\[0.2cm]
\qquad (1) if $ (\Pi_\lambda^h g_n)_K \geq 0$, then $P{\{}\lambda_K, (\Pi_\lambda^h g_n)_K{\}}=0$,\\[0.2cm]
\qquad (2) if $ (\Pi_\lambda^h g_n)_K < 0$, then $P{\{}\lambda_K, (\Pi_\lambda^h g_n)_K{\}}= (\Pi_\lambda^h g_n)_K$,
\item $\lambda_K < 0$, \\[0.2cm]
\qquad (3) $P{\{}\lambda_K, (\Pi_\lambda^h g_n)_K{\}}= (\Pi_\lambda^h g_n)_K$.
\end{itemize}
The optimality conditions are then written as $P{\{}\lambda_K, (\Pi_\lambda^h g_n)_K{\}}=0$. So in the case $(1)$, the constraints are inactive and in the case $(2)$ and $(3)$, the constraints are active.

%\textbf{\huge A revoir}
%
%P. Alart gives a global convergence on a problem between a membrane and an obstacle for a generalised Newton's method on a non-symmetric version of the optimality system of the augmented Lagrangian method \cite{al-cu1988,bica-kozia-08}. The active-set strategy can be re-interpreted as a generalised Newton's method for Alart-Curnier augmented Lagrangian method. In \cite{stad-thesis-04,kunish-stadler-05}, they are showed that the scalability of Newton's method cannot be guarantied. \\

\subsection*{Discrete problem}
\label{subsec: discrete problem}

\noindent The space $V^h$ is spanned by mapped NURBS of type $\hat{N}^p_{\bfi} (\bfz)\circ \varphi_{0,\Gamma_C}^{-1}$ for $\bfi$ belonging to a suitable set of indices. In order to simplify and reduce our notation, we call $A$ as the running index{, of control points associated with the surface $\dis \Gamma_C$,} $A=0 \ldots \caA$ on this basis and set: 
\begin{eqnarray}\label{def:space_disp}
\dis V^h = Span \{ N_A(x), \quad A=0 \ldots \caA \} \cap V.
\end{eqnarray}
Now, we express quantities on the contact interface $\Gamma_C$ as follows: $$\dis \restrr{u}{\Gamma_C} = \sum_{A=1}^{\caA} u_A N_A ,  \qquad  \restrr{\delta u}{\Gamma_C} = \sum_{A=1}^{\caA}  \delta u_A N_A \qquad \textrm{and} \qquad x = \sum_{A=1}^{\caA}  x_A N_A,$$
where %$N_A$ the NURBS basis function include in $N^p(\bfX)$ corresponding to control point $A$, whereas 
$C_A$, $u_A$, $\delta u_A$ and $x_A = \varphi(X_A)$ are the related reference coordinate, displacement, displacement variation and current coordinate vectors. \\

By substituting the interpolations, the normal gap becomes: 
$$g_n = \left[\sum_{A=1}^{\caA} C_A N_A(\zeta)  + \sum_{A=1}^{\caA}  u_A N_A(\zeta) \right] \cdot n .$$
In the previous equation, $\zeta$ are the parametric coordinates of the generic point on $\Gamma_C$. To simplify, we denote for the next of the purpose $\caD g_n[\delta u] =  \delta g_n$. The virtual variation follows as $$\delta g_n = \left[\sum_{A=1}^{\caA}  \delta u_A N_A(\zeta) \right] \cdot n .  $$
In order to formulate the problem in matrix form, the following vectors are introduced: 
$$ \delta \bfu = \begin{bmatrix} \delta u_1\\ \vdots \\ \delta u_{\caA} \end{bmatrix},  \qquad \Delta \bfu = \begin{bmatrix} \Delta u_1\\ \vdots \\ \Delta u_{\caA} \end{bmatrix}, \qquad \bfN = \begin{bmatrix}  N_1(\zeta) n \\ \vdots \\  N_{\caA}(\zeta) n \end{bmatrix}.$$
With the above notations, the virtual variation and the linearized increments can be written in matrix form as follow:
$$ \delta g_n = \delta \bfu^T \bfN , \qquad  \Delta g_n = \bfN^T \Delta \bfu .$$
The contact contribution of the virtual work is expressed as follows: 
\begin{eqnarray}
\begin{array}{l}
\dis \delta W_c = \int_{\Gamma_C} \lambda \delta g_n \ \textrm{d}\Gamma +  \int_{\Gamma_C} \delta \lambda g_n \ \textrm{d}\Gamma .
\end{array}\nonumber
\end{eqnarray}
The discretized contact contribution can be expressed as follows:
\begin{eqnarray}
\begin{array}{ll}
\dis \delta W_c & \dis  = \int_{\Gamma_C}  \sum_{K=1}^{\caK}  \lambda_{K}  {B}_{K} \delta g_n \ \textrm{d}\Gamma +  \int_{\Gamma_C} \sum_{K=1}^{\caK}  \delta \lambda_{K} {B}_{K} g_n \ \textrm{d}\Gamma ,\\[0.4cm] 
& \dis = \sum_{K}  \lambda_{K}   \int_{\Gamma_C} {B}_{K} \delta g_n \ \textrm{d}\Gamma  +   \delta \lambda_{K} \int_{\Gamma_C}  {B}_{K}  g_n \ \textrm{d}\Gamma ,\\[0.4cm] 
& \dis = \sum_{K}  \lambda_{K}   \int_{\Gamma_C} B_K \delta g_n \ \textrm{d}\Gamma  +   \delta \lambda_{K} \int_{\Gamma_C}  B_K  g_n \ \textrm{d}\Gamma ,\\[0.4cm] 
& \dis = \sum_{K} \Big{(}  \lambda_{K}  (\Pi_\lambda^h \delta g_{n})_K  +   \delta \lambda_{K} (\Pi_\lambda^h g_{n})_K \Big{)} K_K ,\\ \end{array}\nonumber
\end{eqnarray}
where $\dis K_K = \int_{\Gamma_C}  B_K \ \textrm{d}\Gamma$.% and the corresponding "local" gaps \eqref{def:l2proj_K} we introduce the following definition of the control point normal gap as the weighted average of the corresponding "local" gaps, with the basis functions as weights: $$\dis (\Pi_\lambda^h  g_{n})_K = \frac{\int_{\Gamma_C} B_K g_n \ \textrm{d}\Gamma}{\int_{\Gamma_C} B_K \ \textrm{d}\Gamma},$$
%and the virtual variation: $$\dis (\Pi_\lambda^h \delta g_{n})_K = \frac{\int_{\Gamma_C} B_K \delta g_n \ \textrm{d}\Gamma}{\int_{\Gamma_C} B_K \ \textrm{d}\Gamma}.$$

\noindent Indeed, we need to resolve a variational inequality. {Using the contact status, we distinguish between the constraints on the control point $K$ are actives, \textit{i.e.} when the contact occurs, and the constraints on the control point $K$ are inactives, \textit{i.e.} when we loose the contact.} \\

\noindent Using active-set strategy on the local gap $ (\Pi_\lambda^h g_{n})_K$ and $\lambda_{K}$, it holds:
\begin{eqnarray}
\begin{array}{l}
\dis \delta W_c = \sum_{K,act} \Big{(}  \lambda_{K} (\Pi_\lambda^h \delta g_{n})_K  +   \delta \lambda_{K} (\Pi_\lambda^h \delta g_{n})_K \Big{)}K_K .
\end{array}\nonumber
\end{eqnarray}
\noindent At the discrete level we proceed as follows:
\begin{itemize}
\item We have $\dis \sum_{K,inact} \delta \lambda_{K} (\Pi_\lambda^h g_{n})_K  \leq 0, \ \forall \delta \lambda_{K}$, \textit{i.e.} $\dis (\Pi_\lambda^h g_{n})_K  \geq 0$ \textit{a.e.} on inactive part.
\item On the active part, it holds $\dis \sum_{K,act} \delta \lambda_{K}  (\Pi_\lambda^h  g_{n})_K  = 0, \ \forall \delta \lambda_{K}$, \textit{i.e.} $\dis (\Pi_\lambda^h g_{n})_K = 0$ \textit{a.e.}.
\item We impose too, $\dis \sum_{K,inact} \lambda_{K} (\Pi_\lambda^h \delta g_{n})_K = 0, \ \forall (\Pi_\lambda^h \delta g_{n})_K$, \textit{i.e.} $\lambda_{K}=0$ \textit{a.e.} on inactive boundary.
\end{itemize}

\noindent For the further developments it is convenient to define the vector of the virtual variations and linearizations for the Lagrange multipliers: $$ \delta \bflambda = \begin{bmatrix} \delta \lambda_{1}\\ \vdots \\ \delta \lambda_{\caK} \end{bmatrix},  \qquad \Delta \bflambda = \begin{bmatrix} \Delta \lambda_{1}\\ \vdots \\ \Delta \lambda_{\caK} \end{bmatrix}, \qquad \bfN_{\lambda,g} = \begin{bmatrix}  (\Pi_\lambda^h g_{n})_{1,act} K_{1,act}  \\ \vdots \\  (\Pi_\lambda^h g_{n})_{\caK,act} K_{\caA, act}  \end{bmatrix}, \qquad \bfB_{\lambda} = \begin{bmatrix}  B_{1}(\zeta)  \\ \vdots \\  B_{\caK}(\zeta)   \end{bmatrix}   .$$

\noindent In the matrix form, it holds: $$\dis \delta W_c = \delta \bfu^T \int_{\Gamma_C} \Big{(} \sum_{K,act} B_K \lambda_{K} \Big{)} \bfN \ \textrm{d}\Gamma + \delta \bflambda^T \bfN_{\lambda,g} , $$
and the residual for Newton-Raphson iterative scheme is obtained as: $$\dis R = \begin{bmatrix} R_u \\ R_\lambda \end{bmatrix}= \begin{bmatrix} \int_{\Gamma_C} \left( \sum_{K,act} B_K \lambda_{K} \right) \bfN \ \textrm{d}\Gamma \\ \bfN_{\lambda,g}  \end{bmatrix} .$$

\noindent The linearization yields: 
\begin{eqnarray}
\begin{array}{l}
\dis \Delta\delta W_c = \int_{\Gamma_C} \Delta \lambda \delta g_n \ \textrm{d}\Gamma +  \int_{\Gamma_C} \delta \lambda \Delta g_n \ \textrm{d}\Gamma .
\end{array}\nonumber
\end{eqnarray}
The active-set strategy and the discretised of contact contribution can be expressed as follows:
\begin{eqnarray}
\begin{array}{ll}
\dis \Delta\delta W_c &\dis= \sum_{K,act}  \sum_{A} \int_{\Gamma_C} \Delta \lambda_{K} B_K N_A \delta u_A \cdot n \ \textrm{d}\Gamma + \int_{\Gamma_C} \delta \lambda_{K} B_K N_A \Delta u_A \cdot n \ \textrm{d}\Gamma ,\\
 &\dis= \delta \bfu^T  \int_{\Gamma_C, act} \bfN \bfB_{\lambda}^T \ \textrm{d}\Gamma \Delta \bflambda + \delta \bflambda^T  \int_{\Gamma_C, act} \bfB_{\lambda} \bfN^T \ \textrm{d}\Gamma \Delta \bfu.
\end{array}\nonumber
\end{eqnarray}

%%%%%%%%%%%%%%%%%%%%%%%%%%%%%%%%%%%%%%%%%%%%%%%%
 %Add in table of contents
% \addcontentsline{toc}{section}{Appendices}

%\newpage

%\section*{{Appendix 1. Large deformation}}
%\label{sec:ALM}
%\input{large_def.tex}
%%%%%%%%%%%%%%%%%%%%%%%%%%%%%%%%%%%%%%%%%%%%%%%%

\newpage

%%%%%%%%%%%%%%%%%%%%%%%%%%%%%%%%%%%%%%%%%%%%%%%%
 \addcontentsline{toc}{section}{Bibliography}
  
\bibliographystyle{siam}

%%%%%%%%%%%%%%%%%%%%%%%%%%%%%%%%%%%%%%%%%%%%%%%%

\end{document}